\documentclass{amsart}
\usepackage{amsmath}
\usepackage{amssymb} 
\usepackage{amscd} 
\usepackage{graphicx} 
\usepackage{epsfig}
\usepackage[dvips]{color}
\usepackage{bm}

\newtheorem{theorem}{Theorem}[section]
\newtheorem{lemma}[theorem]{Lemma}

\newtheorem{proposition}[theorem]{Proposition}

\newtheorem{claim}[theorem]{Claim}

\theoremstyle{definition}
\newtheorem{definition}[theorem]{Definition}

\theoremstyle{remark}
\newtheorem{remark}[theorem]{Remark}

\numberwithin{equation}{section}
\numberwithin{figure}{section}
\numberwithin{table}{section}

\begin{document}

\title[Networking Seifert Surgeries on Knots IV]
{Networking Seifert Surgeries on Knots IV:\\
Seiferters and branched coverings}

\author[A. Deruelle]{Arnaud Deruelle}
\address{Institute of Natural Sciences, 
Nihon University, 
Tokyo 156--8550, Japan}
\email{aderuelle@math.chs.nihon-u.ac.jp}

\author[M. Eudave-Mu\~noz]{Mario Eudave-Mu\~noz}
\address{Instituto de Matematicas, Universidad Nacional Aut\'onoma de M\'exico, Circuito Exterior, Ciudad Universitaria 04510 M\'exico DF, Mexico, and CIMAT, Guanajuato, Mexico}
\email{mario@matem.unam.mx}
\thanks{The second author was partially supported by PAPIIT-UNAM grant IN102808. }

\author[K. Miyazaki]{Katura Miyazaki}
\address{Faculty of Engineering, Tokyo Denki University, Tokyo 101--8457, 
Japan}
\email{miyazaki@cck.dendai.ac.jp}

\author[K. Motegi]{Kimihiko Motegi}
\address{Department of Mathematics, Nihon University, 
Tokyo 156--8550, Japan}
\email{motegi@math.chs.nihon-u.ac.jp}
\thanks{
The last author has been partially supported by JSPS Grants--in--Aid for Scientific 
Research (C) (No.17540097 and No.21540098), The Ministry of Education, Culture, Sports, Science and Technology, Japan and Joint Research Grant of Institute of Natural Sciences at 
Nihon University for 2011. }

\date{}


\begin{abstract}
A Seifert surgery is an integral surgery on a knot in
$S^3$ producing a Seifert fiber space $M$
which may contain an exceptional fiber of index $0$.
The Seifert Surgery Network is a $1$--dimensional complex
whose vertices correspond to Seifert surgeries;
its edges correspond to single twistings along 
``seiferters" or ``annular pairs of seiferters".  
One problem of the network is
whether there is a path from each vertex
to a vertex on a torus knot,
the most basic Seifert surgery.
We give a method to find seiferters and annular pairs of seiferters for 
Seifert surgeries obtained
by taking two--fold branched covers of tangles.
Concerning three infinite families of Seifert surgeries
obtained by the second author via branched covers,
we find explicit paths in the network
from such surgeries
to Seifert surgeries on torus knots.
\end{abstract}

\maketitle

\section{Introduction}
\label{section:Introduction}

Let $K$ be a knot in the $3$--sphere $S^3$ and $m$ an integer. 
If the result $K(m)$ of $m$--Dehn surgery on $K$
is a Seifert fiber space which may have a fiber of index zero, 
then we call the pair $(K, m)$ a \textit{Seifert surgery}.  
As shown in  \cite[Proposition~2.8]{DMM1} if $K(m)$ admits a Seifert fibration with 
fiber of index zero, then it is a lens space or a connected sum of two lens spaces. 
For any nontrivial torus knot $T_{p, q}$, 
the Seifert surgery $(T_{p, q}, pq)$ is such an example; 
$T_{p, q}(pq) \cong L(p, q) \sharp L(q, p)$.  

\begin{figure}[!ht]
\begin{center}
\includegraphics[width=0.4\linewidth]{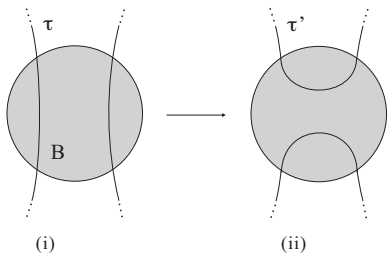}
\caption{Tangle replacement.}
\label{inftytangle}
\end{center}
\end{figure}

Many examples of Seifert surgeries are constructed by using the Montesinos trick 
(\cite{BH}, \cite{EM1,EM2}, \cite{BZ1}).  
Let $\tau$ be a trivial knot in $S^3$,
and $B$ a $3$--ball such that 
$(B, B \cap \tau)$ is a $2$--string trivial tangle
as in Figure~\ref{inftytangle}(i). 
Suppose that $\tau$ is changed to a Montesinos link 
or a Montesinos-m link $\tau'$ after replacing 
$(B, B \cap \tau)$ with another trivial tangle;
for example, see Figure~\ref{inftytangle}(ii).  

Since $\tau$ is a trivial knot
and $(B, B\cap \tau)$ is a $2$--string trivial tangle, 
the two--fold branched cover of $S^3$ along $\tau$ is $S^3$
and the preimage $\widetilde{B}$ of $B$ is a solid torus; 
denote the core of the solid torus by $K$.
Replacing the tangle $(B, B\cap \tau)$ downstairs
corresponds to
removing the solid torus $\widetilde{B} =N(K)$
and re-attaching upstairs,
i.e.\ Dehn surgery on $K$.
This observation is referred to as the Montesinos trick \cite{Mon}. 
In our case,
since $\tau'$ is a Montesinos link
(resp.\ Montesinos--m link),
the two--fold branched cover of $S^3$ along $\tau'$
is a Seifert fiber space over $S^2$ 
(resp.\ $\mathbb{R}P^2$)\cite{Montesinos}. 
It follows that $K$ has a Seifert surgery.
For details, see Section~\ref{section:covering}. 

In \cite{EM2},
the second author solved tangle equations,  
and found three infinite families of Seifert surgeries
by using branched covers. 
We denote these families by $\mathcal{EM}\mathrm{I},
\mathcal{EM}\mathrm{II}, \mathcal{EM}\mathrm{III}$. 
In this paper, we study these surgeries
from a viewpoint of the Seifert Surgery Network
introduced in \cite{DMM1}. 
In \cite{DMM1} we define relationships among Seifert surgeries,
and draw a global picture of Seifert surgeries. 
To do this
we have introduced seiferters and the Seifert Surgery Network, 
a $1$--dimensional complex whose vertices correspond to Seifert surgeries. 

\begin{definition}[seiferter]
A knot $c$ in $S^3 - N(K)$ is called a \textit{seiferter} for 
a Seifert surgery $(K, m)$
if $c$ enjoys the following two properties.
\begin{enumerate}
\item $c$ is unknotted in $S^3$.
\item $c$ becomes a fiber in
a Seifert fibration of $K(m)$.
\end{enumerate}
\end{definition}

Let $c_1, c_2$ be seiferters for $(K, m)$
such that $c_1$ and $c_2$ become fibers in some
Seifert fibration of $K(m)$ simultaneously.
Then we call $\{c_1, c_2\}$
\textit{a pair of seiferters} for $(K, m)$.
Furthermore,
if $c_1, c_2$ cobound an annulus in $S^3$,
then we call $\{c_1, c_2\}$
\textit{an annular pair of seiferters}.

\begin{remark}
\label{(c_1, c_2)}
Let $A$ be an annulus cobounded by $c_1, c_2$.
In \cite{DMM1},
an annular pair $\{c_1, c_2\}$ is defined
to be an ordered pair of $c_1$ and $c_2$ to specify
the direction of twist along the annulus $A$.
In this paper,
annular pairs are presented
as ordered pairs only when we perform twisting
along annuli.
If $\mathrm{lk}(c_1, c_2) =l$
where $c_i$ are oriented so as to be
homologous in $A$,
then we define \textit{$p$--twist along $(c_1, c_2)$}
to be performing $(-\frac{1}{p}+l)$--surgery on $c_1$
and $(\frac{1}{p}+l)$--surgery on $c_2$.
This pair of surgeries is equivalent
to twisting $p$ times along the annulus $A$. 
See \cite[Definition~2.32]{DMM1} for details.
\end{remark}

For a Seifert surgery $(K, m)$ with a seiferter $c$ 
(resp.\ an annular pair $\{ c_1, c_2 \}$),
let $K_p$ and $m_p$ be the images of $K$ and $m$
under $p$--twist along $c$ (resp.\ $( c_1, c_2 )$),
respectively. 
The key fact is that $(K_p, m_p)$ remains a Seifert surgery for any integer $p$, 
and $c$ (resp.\ $\{ c_1,\ c_2 \}$) remains a seiferter 
(resp.\ an annular pair) for $(K_p, m_p)$ \cite[Propositions~2.6 and 2.33]{DMM1}. 
We say that $(K_p, m_p)$ is obtained from $(K, m)$ by twisting along 
$c$ (resp.\ $\{ c_1, c_2 \}$). 
See Diagram 1 below. 

\begin{eqnarray*}
\begin{CD}
	(K, m) @>\textrm{twist along } c\ (\textrm{resp.\ } (c_1, c_2)) >> (K_p, m_p) \\
@V{m \textrm{--surgery on } K}VV
		@VV{m_p \textrm{--surgery on } K_p}V \\
	K(m) @>>\textrm{surgery on } c\ (\textrm{resp.\ } (c_1, c_2)) > K_p(m_p)
\end{CD}
\end{eqnarray*}
\vskip 0.2cm
\begin{center}
\textsc{Diagram 1.} 
\end{center}

\medskip

In the Seifert Surgery Network,
two vertices (i.e. two Seifert surgeries) are connected by an edge if
one is obtained from the other by a single twist along a seiferter
or an annular pair of seiferters \cite[Subsection~2.4]{DMM1}. 
For a Seifert surgery $(K, m)$ with a seiferter $c$
(resp.\ an annular pair of seiferters $\{ c_1, c_2 \}$),
twisting $(K, m)$ successively along $c$ (resp.\ $(c_1, c_2)$)
naturally generates a 1--dimensional subcomplex. 

The purpose of this paper is to prove:

\begin{theorem}
\label{connected}
For each Seifert surgery in $\mathcal{EM}\mathrm{I}$, 
$\mathcal{EM}\mathrm{II}$, and $\mathcal{EM}\mathrm{III}$,
there is a path to a Seifert surgery on a torus knot 
in the Seifert Surgery Network.
\end{theorem}

Since the exteriors of torus knots are Seifert fibered,
Seifert surgeries on torus knots are the most basic Seifert surgeries
and well understood.  
Hence,
if we trace a given Seifert surgery $(K, m)$ back to
a Seifert surgery on a torus knot, 
then we can regard the surgery on the torus knot
as an origin of $(K, m)$ in the network. 
In \cite{DMM1}, 
we found out origins of several Seifert surgeries, 
and in \cite{DMM2} we gave explicit paths
from the Berge's lens surgeries \cite{Berge2} to 
Seifert surgeries on torus knots.

In Section~\ref{section:seiferters and tangles}, 
we give a general method
to find seiferters and annular pairs of seiferters for  
Seifert surgeries obtained by using branched coverings 
(Theorems~\ref{seiferter for covering knots} and \ref{annular pair for covering knots}). 
In Section~\ref{section:examples},
we apply results in
Section~\ref{section:seiferters and tangles}
to find explicit paths from such Seifert surgeries to surgeries on torus knots.
In Subsections~\ref{subsection:Klnp},
\ref{subsection:Klmnp}, \ref{subsection:KABCi},
Theorem~\ref{connected} is proved
for $\mathcal{EM}\mathrm{I}$,
$\mathcal{EM}\mathrm{II}$,
$\mathcal{EM}\mathrm{III}$, respectively.

\bigskip

\section{Tangles, branched coverings and Seifert surgeries}
\label{section:covering}

A \textit{tangle} $(B, t)$ is a pair of a $3$--ball $B$ and two disjoint arcs $t$ properly embedded in $B$. 
A tangle $(B, t)$ is trivial if 
there is a pairwise
homeomorphism from $(B, t)$ to $(D^2 \times I, \{x_1, x_2\}\times I)$, where $x_1, x_2$ are distinct points.
For tangles $(B, t)$ and $(B, t')$ with $\partial t = \partial t'$,
we say that they are \textit{equivalent} 
if there is a pairwise homeomorphism 
$h : (B, t) \to (B, t')$
satisfying $h|_{\partial B} =$ id. 

\begin{figure}[htbp]
\begin{center}
\includegraphics[width=0.9\linewidth]{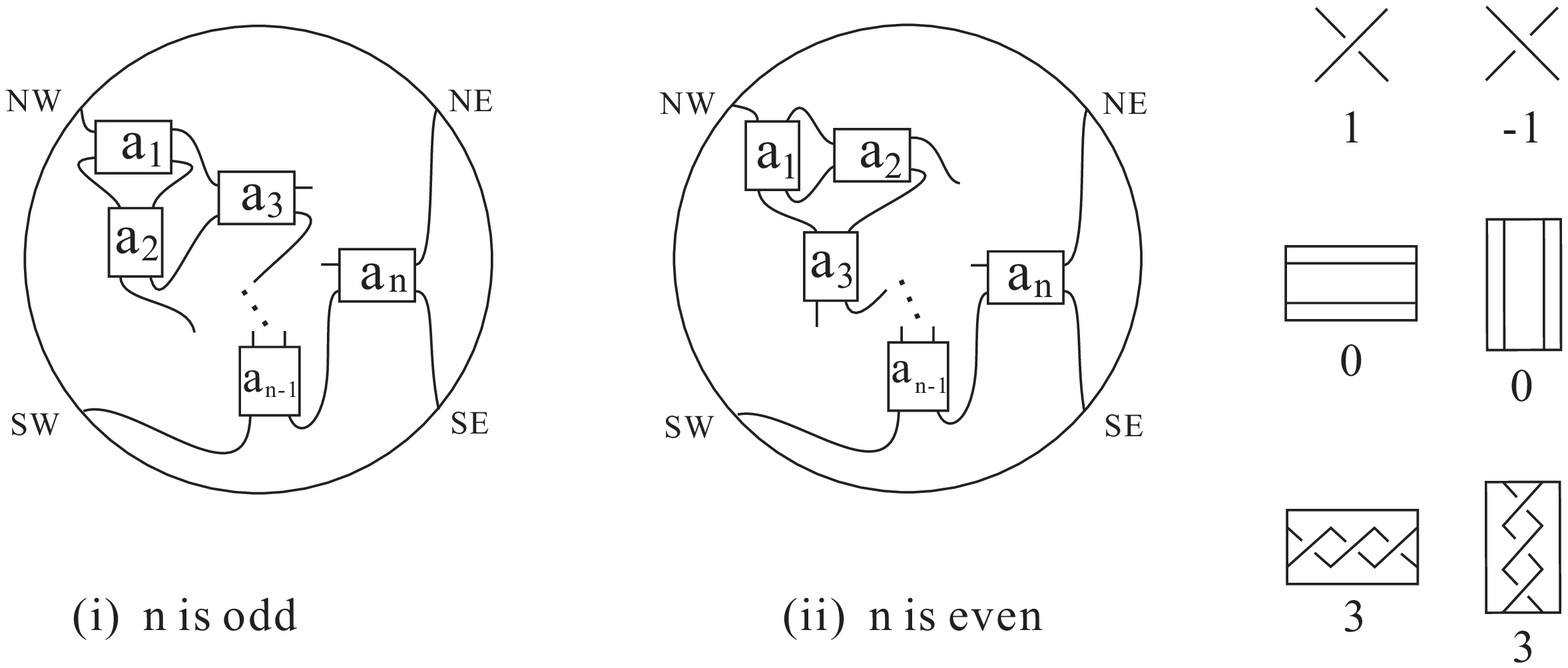}
\caption{Rational tangles.}
\label{rtangle}
\end{center}
\end{figure}

Let $U$ be the unit 3-ball in $\mathbb{R}^3$,
and take 4 points NW, NE, SE, SW on the boundary of $U$
so that
$\mathrm{NW} = (0,-\alpha, \alpha),
\mathrm{NE} = (0, \alpha, \alpha),
\mathrm{SE} = (0, -\alpha, \alpha),
\mathrm{SW} = (0, -\alpha, -\alpha)$,
where $\alpha =\frac{1}{\sqrt{2}}$.
A tangle $(U, t)$ is a \textit{rational tangle}
if it is a trivial tangle with $\partial t =
\{\mathrm{NW, NE, SE, SW}\}$.
We can construct rational tangles from sequences of integers 
$a_1, a_2, \dots, a_n$ as shown in Figure~\ref{rtangle}, 
where the last horizontal twist $a_n$ may be $0$.
We consider that the tangle diagrams in Figure~\ref{rtangle}
is drawn on the $yz$--plane. 
Denote by $R(a_1, a_2, \dots, a_n)$ the associated rational tangle.  

Each rational tangle can be parametrized by 
$r \in \mathbb{Q} \cup \{ \infty\}$, 
where the rational number $r$ is given by
the continued fraction below. 
Thus we denote the rational tangle corresponding to $r$ by $R(r)$. 

$$ r\ =\ a_{n} + \cfrac 1 {a_{n-1}+ \cfrac 1 {
               \begin{array}{clr}
               \ & & \\[-5pt]
               \hspace{-25pt} \ddots & & \\[-10pt]
                      & \ \ \hspace{-5pt} a_2+ \cfrac{1}{a_1}
               \end{array}
               }}$$
               
\begin{figure}[!ht]
\begin{center}
\includegraphics[width=0.35\linewidth]{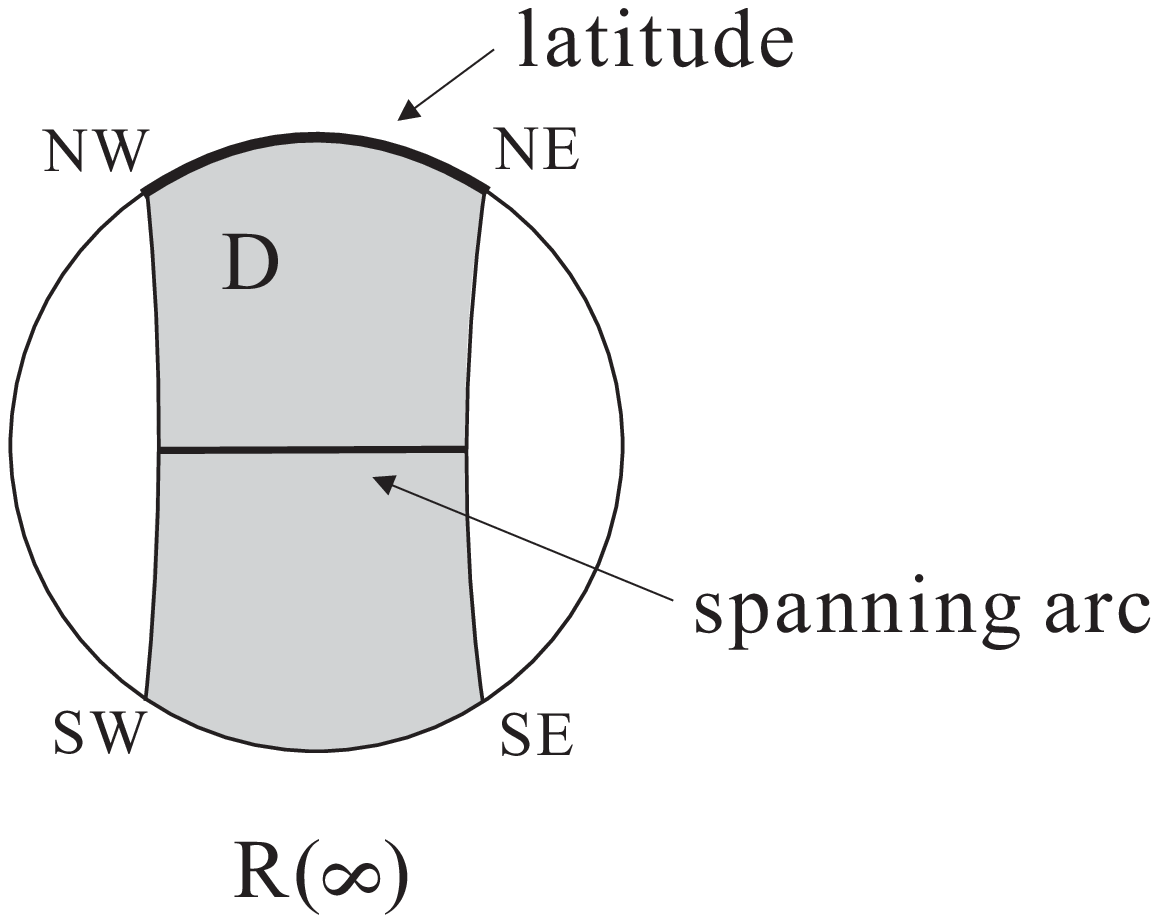}
\caption{}
\label{spanningarc}
\end{center}
\end{figure}

Let $(U, t)$ be the rational tangle $R(\infty)$.
Considering $t$ is embedded in the $yz$--plane,
take the disk $D$ in the $yz$--plane such that
$\partial D$ is the union of $t$ and 
two arcs in $\partial U$: one connects
NW and NE, and the other connects SW and SE.
We call an arc in $D$ connecting the components of the interior of $t$
a \textit{spanning arc},
and the arc $D\cap \partial U$ connecting
NW and NE \textit{the latitude of $R(\infty)$}. 
See Figure~\ref{spanningarc}. 
The two--fold cover $\widetilde{U}$
of $U$ branched along $t$ is a solid torus.
Note that the preimages of the spanning arc
and the latitude are a core and a longitude $\lambda$
of the solid torus, respectively.
A meridian of a rational tangle $R(r)=(U, t')$ is
a simple closed curve in $\partial U -t'$ which bounds a disk
in $U -t'$ and a disk in $\partial U$ meeting $t'$
in two points.
Let $\mu_r (\subset \partial \widetilde{U})$
be a lift of a meridian of $R(r)$;
then $\mu_{\infty}$ is a meridian of the solid torus $\widetilde{U}$.
Furthermore, we note the following well-known fact.

\begin{lemma}
\label{-p/q}
Under adequate orientations we have
$[\mu_r] = -p[\mu_{\infty}] +q[\lambda]
\in H_1(\partial \widetilde{U})$,
where $r = \frac{p}{q}$ and $[\mu_{\infty}]\cdot[\lambda] =1$.
\end{lemma}

Let $(B, t)$ be a tangle such that
$B \subset S^3 ( = \mathbb{R}^3\cup\{\infty\} )$ is
the complement of the unit 3-ball $U$,
and $\partial t = \{ \mathrm{NW, NE, SE, SW} \}$.
We denote by $(B, t) + R(r)$ 
the knot or link in $S^3$ 
formed by the union of the strings of the tangles,
and let $\pi_r: X_r \to S^3 =B\cup U$ be
the two--fold cover branched along $(B, t) + R(r)$.
We say that $(B, t)$ is \textit{trivializable}
if  $(B, t) + R(\infty)$ is a trivial knot in $S^3$. 
If $(B, t) + R(r)$ is a trivial knot for some $r\in \mathbb{Q}$, 
then an ambient isotopy of $B$ changes
$(B, t)$ to a trivializable tangle.

Suppose that $(B, t) + R(\infty)$ is a trivial knot. 
Then the two--fold branched cover $X_{\infty}$ is the 3--sphere,
and the preimage of the spanning arc $\kappa$ for $R(\infty)$
is a knot in $X_{\infty} =S^3$,
which we call the \textit{covering knot} of $(B, t)$.
The exterior of the covering knot $K$
is $\pi_{\infty}^{-1}(B)$.
For $(B, t) +R(\infty)$
a replacement of $R(\infty)$ by a rational tangle $R(s)$
is called \textit{$s$--untangle surgery}
on $(B, t) +R(\infty)$. 
Performing untangle surgery downstairs
corresponds to replacing the solid torus $\pi_{\infty}^{-1}(U)$
by $\pi_s^{-1}(U)$ upstairs,
i.e.\ Dehn surgery on the covering knot $K$.
We denote the surgery slope by $\gamma_s$;
it is represented by a lift
of a meridian of $R(s)$. 
We say that $\gamma_s$ is the \textit{covering slope} of $s$. 
See the commutative diagram below. 

\begin{eqnarray*}
\begin{CD}
	S^3 @>\gamma_s\textrm{--surgery on } K >> K(\gamma) \\
@V{ \textrm{two--fold branched cover}}VV
		@VV{\textrm{two--fold branched cover}}V \\
	(B, t)\cup R(\infty) @>>s\textrm{--untangle surgery} > (B, t)\cup R(s)
\end{CD}
\end{eqnarray*}
\vskip 0.2cm
\begin{center}
\textsc{Diagram 2.} Montesinos trick 
\end{center}

\begin{remark}
\label{preferred framing}
If the preimage of the latitude of $R(\infty)$
is a preferred longitude of the covering knot $K$,
then by Lemma~\ref{-p/q}
the covering slope $\gamma_s$,
where $s = \frac{p}{q}$,
is $-\frac{p}{q}$ in terms of
a preferred meridian--longitude pair of $K$.
\end{remark}

For a link $L$ and an arc $\tau$ with $\tau \cap L = \partial \tau$
we perform an untangle surgery along $\tau$ as follows.
First take a regular neighborhood $N(\tau)$ of $\tau$ so that
$(N(\tau), N(\tau) \cap L)$ is a trivial tangle.
Then, identifying the trivial tangle $T =(N(\tau), N(\tau) \cap L)$
with the rational tangle $R(\infty)$,
we can replace $R(\infty) = T$ by a rational tangle $R(s)$;
this operation is called $s$--untangle surgery of $L$ along $\tau$.
Note that the definition of $s$-untangle surgery along $\tau$
relies on the identification of $T$ with $R(\infty)$.
If $L$ is a trivial knot, 
the two--fold branched cover of $S^3$ along $L$ is $S^3$,
and the preimage of $\tau$ is a knot, which we call the covering knot
of $\tau$.
Then, as before,
performing $s$--untangle surgery along $\tau$ downstairs
corresponds to performing Dehn surgery on
the covering knot upstairs;
we call its surgery slope the covering slope. 
For two disjoint arcs $\tau_1, \tau_2$ with $\tau_i \cap L = \partial \tau_i$, a pair of $s_i$--untangle surgeries along $\tau_i$
is called \textit{
$(s_1, s_2)$--untangle surgery along $(\tau_1, \tau_2)$}.

\begin{figure}[!ht]
\begin{center}
\includegraphics[width=0.5\linewidth]{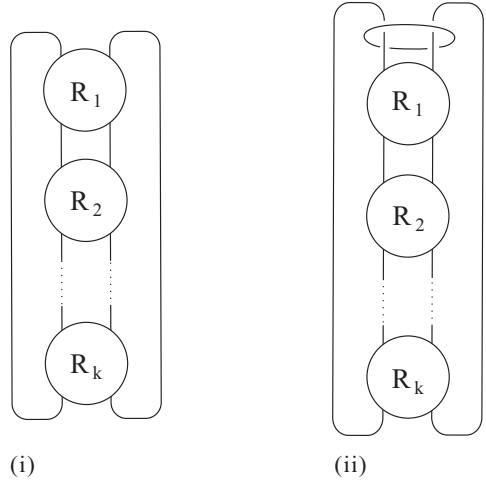}
\caption{Standard positions of Montesinos link(--m) links.}
\label{Montesinos}
\end{center}
\end{figure}

A \textit{Montesinos link} $M(R_1, \cdots, R_k)$
(resp.\ \textit{Montesinos--m link} $Mm(R_1, \cdots, R_k)$)
is a link which has a diagram in Figure~\ref{Montesinos}(i)
(resp.\ (ii)),
where $R_i$ are rational tangles as shown in Figure~\ref{rtangle}.
We call the diagrams in Figure~\ref{Montesinos} 
\textit{standard positions} of Montesinos(--m) links. 
If $R_i$ corresponds to $r_i\in \mathbb{Q}\cup \{\infty\}$ for $i=1, \ldots, k$,
then we often write
$M(r_1, \ldots, r_k)$ or $Mm(r_1, \ldots, r_k)$.
Let $X$ be the two--fold branched cover of $S^3$ along
a Montesinos link $M(r_1, \ldots, r_k)$ (resp.\ 
a Montesinos--m link $Mm(r_1, \ldots, r_k)$).
Then $X$ admits a Seifert fibration in which
the preimage of $B_i$, where $R_i =(B_i, r_i)$,
is a fibered solid torus and its core has
Seifert invariant $-\frac{1}{r_i}$ and index $|p_i|$,
where $r_i = \frac{p_i}{q_i}$.
Hence, $X =S^2(-\frac{1}{r_1}, \ldots, -\frac{1}{r_k})$ 
(resp.\ $\mathbb{R}P^2(-\frac{1}{r_1}, \ldots, -\frac{1}{r_k})$). 
See \cite{Montesinos}.

Let $(B, t)$ be a trivializable tangle such that
$(B, t) + R(s)$ is a Montesinos(--m) link
for some rational number $s$.
The two--fold branched cover $X_s$,
which is a Seifert fiber space as we see above,
is obtained from $S^3$ by
$\gamma_s$--surgery on the covering knot $K$ of $(B, t)$.
In this manner, we obtain a Seifert surgery $(K, \gamma_s)$.

\section{Seiferters, annular pairs of seiferters and tangles}
\label{section:seiferters and tangles}

In this section we assume that
a tangle $(B, t)$ satisfies the following conditions: 
\begin{itemize}
\item
$L_{\infty} = (B, t) + R(\infty)$ is a trivial knot in $S^3$, and 

\item
$L_s = (B, t) + R(s)$ is a Montesinos link $M(R_1, \ldots, R_k)$
or a Montesinos--m link $Mm(R_1, \ldots, R_k)$.
\end{itemize}

As in the previous section,
we denote by $K$
the covering knot of the trivializable tangle $(B, t)$,
and by $\gamma_s$ the covering slope corresponding to
the replacement of $R(\infty)$ with $R(s)$.
We let $\pi_{\infty}: X_{\infty} \to S^3$ (resp.\
$\pi_s: X_s\to S^3$) denote the two--fold cover of $S^3$
branched along $L_{\infty}$ (resp.\ $L_s$).
The Montesinos(--m) link 
$L_s$ can be deformed into a standard position as in Figure~\ref{Montesinos}.
We define a leading arc of a rational tangle.
The preimage of a leading arc then becomes a Seifert fiber in $X_s$.

\begin{figure}[!ht]
\begin{center}
\includegraphics[width=0.7\linewidth]{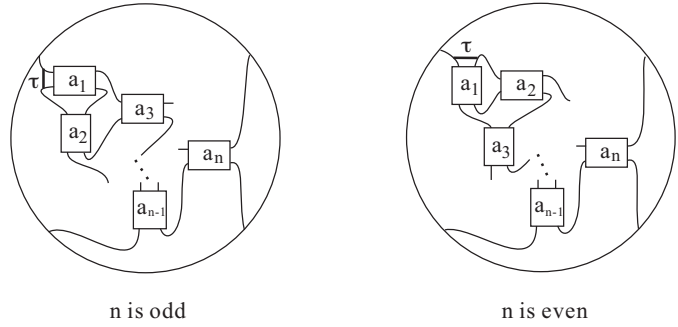}
\caption{Leading arcs in rational tangles.}
\label{fig:leadingarc}
\end{center}
\end{figure}

\begin{definition}[leading arc]
\label{leadingarc}
Let $\tau$ be an arc in a rational tangle $R = R(a_1, \dots, a_n)$ 
as depicted in Figure~\ref{fig:leadingarc}. 
Then we call $\tau$ a {\it leading arc} of $R$.
\end{definition}

\begin{lemma}
\label{fiber}
Let $\tau$ be a leading arc of a rational tangle
$R_i = R(\frac{p_i}{q_i})$ in a standard position of $L_s$. 
Then $c =\pi_s^{-1}(\tau)$ is a fiber of index $|p_i|$
in a Seifert fibration of $X_s$. 
\end{lemma}

\textsc{Proof.}
Since $R_i = (B_i, t_i)$ is a rational tangle,
by an ambient isotopy of $B_i$ we can deform $(B_i, t_i)$ and $\tau$ 
to $R(\infty)$ and a spanning arc for $R(\infty)$
as in Figure~\ref{spanningarc}.
Hence $c$ is a core of the solid torus $\pi_s^{-1}(B_i)$.
This implies the desired result.
\hspace*{\fill} $\square$(Claim~\ref{fiber})

\begin{remark}
\label{indexp} 
In Lemma~\ref{fiber}, 
if $|p_i|=1$, then $c$ is a regular fiber in $X_s$. 
If $p_i = 0$, i.e.\ $R_i = R(0)$, 
then $L_s$ is not a Montesinos(--m) link in the usual sense, 
and $X_s$ is a connected sum of lens spaces; 
$c$ is a degenerate fiber in $X_s$. 
\end{remark}

\begin{theorem}[seiferters for covering knots]
\label{seiferter for covering knots}
Let $\tau$ be an arc in $B$ such that $\tau \cap t = \partial \tau$.  
Assume that after an isotopy of $\tau \cup L_s$, 
$\tau$ is a leading arc of some $R_i$ in a standard position of $L_s = (B, t) + R(s)$.  
Assume that some nontrivial untangle surgery along $\tau$ 
preserves the triviality of $L_{\infty} = (B, t) + R(\infty)$.
Then the following hold. 
\begin{enumerate}
\item
The preimage $c = \pi_{\infty}^{-1}(\tau)$ is a seiferter for the Seifert surgery $(K, \gamma_s)$.

\item The above untangle surgery along $\tau$ corresponds to twisting along the seiferter $c$ in $S^3$.
\end{enumerate}
\end{theorem}

\textsc{Proof.}
By Lemma~\ref{fiber} $c$ is a Seifert fiber in $K(\gamma_s)$. 
To prove (1) it remains to show that $c$ is a trivial knot in $S^3$. 
By the second assumption 
some nontrivial Dehn surgery of $S^3$ on $c$ yields $S^3$.
It then follows from \cite{GL2} that
$c$ is a trivial knot in $S^3$, and thus
the Dehn surgery on $c$ is $\frac{1}{n}$--surgery for some integer $n$.
This proves (1) and (2). 
\hspace*{\fill} $\square$(Theorem~\ref{seiferter for covering knots})

\begin{remark}
\label{framing}
If for some integer $n_0$ $(|n_0| > 2)$, 
$\frac{1}{n_0}$--untangle surgery along $\tau$ preserves the triviality of 
$(B, t) + R(\infty)$, 
then the preimage of
the latitude of $R(\infty)=(N(\tau), N(\tau)\cap L_{\infty})$ 
is a preferred longitude of $c$. 
In fact, 
the covering slope corresponding to $\frac{1}{n_0}$ is
$x - \frac{1}{n_0} = \frac{xn_0 -1}{n_0}$ for some integer $x$.  
The proof of Theorem~\ref{seiferter for covering knots}(1) shows that 
$|xn_0 -1| = 1$. 
Since $|n_0| > 2$, $x$ must be zero. 
Thus for any integer $n$, 
$\frac{1}{n}$--untangle surgery along $\tau$ corresponds to 
$-\frac{1}{n}$--surgery on $c$. 
\end{remark}

\begin{figure}[htbp]
\begin{center}
\includegraphics[width=0.25\linewidth]{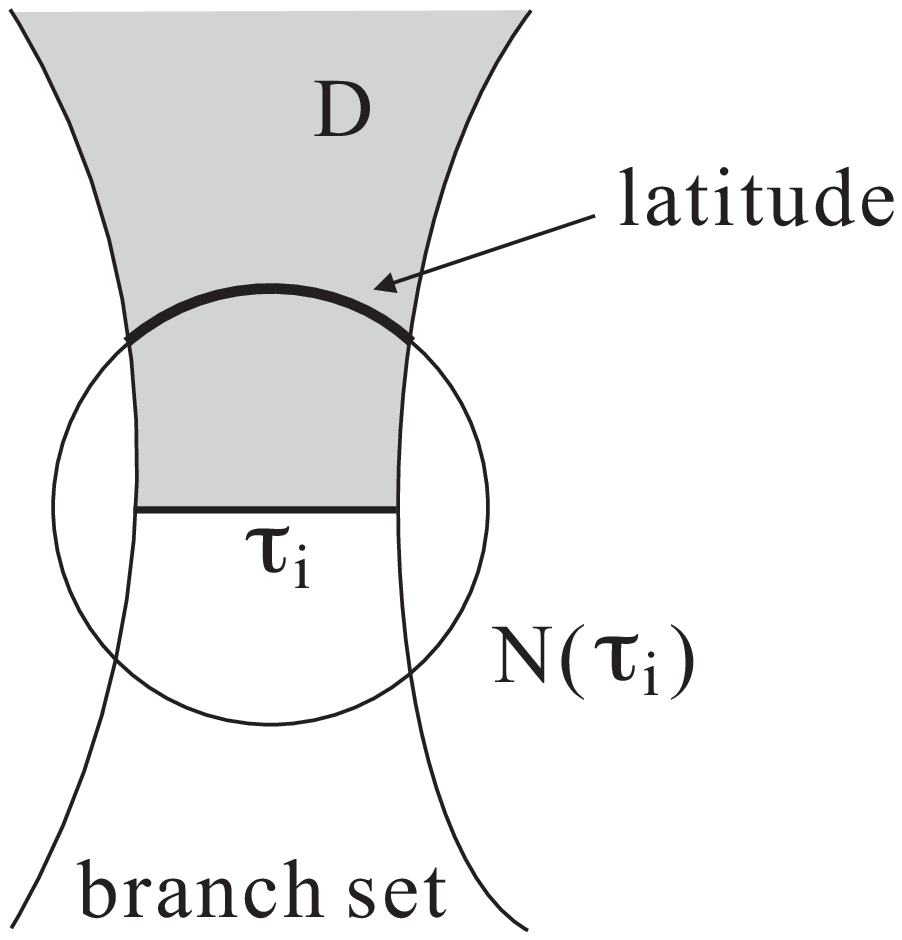}
\caption{}
\label{latitude}
\end{center}
\end{figure}

\begin{theorem}
[annular pairs of seiferters for covering knots]
\label{annular pair for covering knots}
Let $\tau_1$ and $\tau_2$ be disjoint arcs in $B$ such that 
$\tau_i \cap t = \partial \tau_i$.  
Assume that {\upshape (i)} $\tau_1$ and $\tau_2$ become
leading arcs of $R_i$ and $R_j$ $(i\neq j)$ 
respectively
in a standard position of the Montesinos(--m) link $L_s$, 
and {\upshape (ii)} for each $i = 1, 2$
some nontrivial untangle surgery along $\tau_i$
preserves the triviality of $L_{\infty}$. 
Then by Theorem~\ref{seiferter for covering knots},  
the preimages $c_i =\pi_{\infty}^{-1}(\tau_i)$ $(i=1, 2)$ 
form a pair of seiferters for $(K, \gamma_s)$. 
Suppose further that there is a rectangle $D$ such that
$\tau_1$ and $\tau_2$ are two opposite sides of $D$, 
the other sides are contained in $L_{\infty}$, and 
$\mathrm{int}D \cap L_{\infty} = \emptyset$. 
Then the following hold.
\begin{enumerate}
\item
$\{c_1, c_2\}$ is an annular pair of seiferters
for $(K, \gamma_s)$. 
\item
Isotope $\mathrm{int}D$
so that $D\cap N(\tau_i)$ is a disk for each $i$,
as in Figure~\ref{latitude}.
Identify $(N(\tau_i), N(\tau_i)\cap L_{\infty})$ with $R(\infty)$ 
so that the arc $D \cap \partial N(\tau_i)$ becomes
the latitude of $R(\infty)$.
Then, under this identification, 
$(\frac{1}{n},\ \frac{-1}{n})$--untangle surgery along
$(\tau_1, \tau_2)$ 
corresponds to $n$--twist along the annular pair of
seiferters $(c_1, c_2)$.  
\end{enumerate}
\end{theorem}

\textsc{Proof.}
We first show that 
the preimage $\pi_{\infty}^{-1}( D )$ is an annulus
cobounded by $c_1$ and $c_2$.
The rectangle $D$ intersects the branch set $L_{\infty}$
in two opposite sides.
Let us consider the rectangle $D'= D - \mathrm{int}N(L_{\infty})$. 
Then the preimage of $D'$ (under unbranched covering) consists of two disjoint disks. 
Completion along the branch set,
we obtain an annulus cobounded by the seiferters $c_1$ and $c_2$
as claimed in (1).

Let $A$ be the annulus $\pi_{\infty}^{-1}(D)$. 
Let $l$ be the linking number between $c_1$ and $c_2$,
where $c_1$ and $c_2$ are oriented so as to be homologous in $A$.
After the isotopy in assertion~(2)
$A$ intersects the boundary of
$N(c_i) =\pi_{\infty}^{-1}( N(\tau_i ))$
in a longitude $\lambda$ of slope $l$.
By the identification given in assertion~(2)
$\lambda$ is also the preimage of the latitude for $R(\infty)$.
It follows from Lemma~\ref{-p/q} that
the covering slope of $\frac{\pm1}{n}$
is $\mp [\mu]+n[\lambda] \in H_1(\partial N(c_i))$
where $\mu$ is a meridian of $N(c_i)$ with $[\mu] \cdot [\lambda] = 1$.  
Hence, $(\frac{1}{n},\ \frac{-1}{n})$--untangle surgery along
$(\tau_1, \tau_2)$ corresponds to
$(-\frac{1}{n}+l, \frac{1}{n}+l)$--surgery along $(c_1, c_2)$
in terms of preferred meridian-longitude coordinates, 
i.e.\ $n$--twist along $(c_1, c_2)$ by Remark~\ref{(c_1, c_2)}. 
\hspace*{\fill} $\square$(Theorem~\ref{annular pair for covering knots})

\section{Seifert surgeries on covering knots}
\label{section:examples}

Let us turn to specific examples of Seifert surgeries obtained by using the Montesinos trick.  
In \cite{EM1} the second author gave explicit infinite families
of Seifert surgeries on hyperbolic knots using the Montesinos trick.
In \cite{EM2},
by solving tangle equations he expanded these families
and found two more infinite families.

\subsection{The first family of Seifert surgeries 
$\mathcal{EM}\mathrm{I}$}
\label{subsection:Klnp}

The first family $\mathcal{EM}\mathrm{I}$
consists of Seifert surgeries on knots $K(l,n,p)$,
which are the covering knots of the tangles $T(l, n, p)$
below.

\begin{figure}[htbp]
\begin{center}
\includegraphics[width=0.55\linewidth]{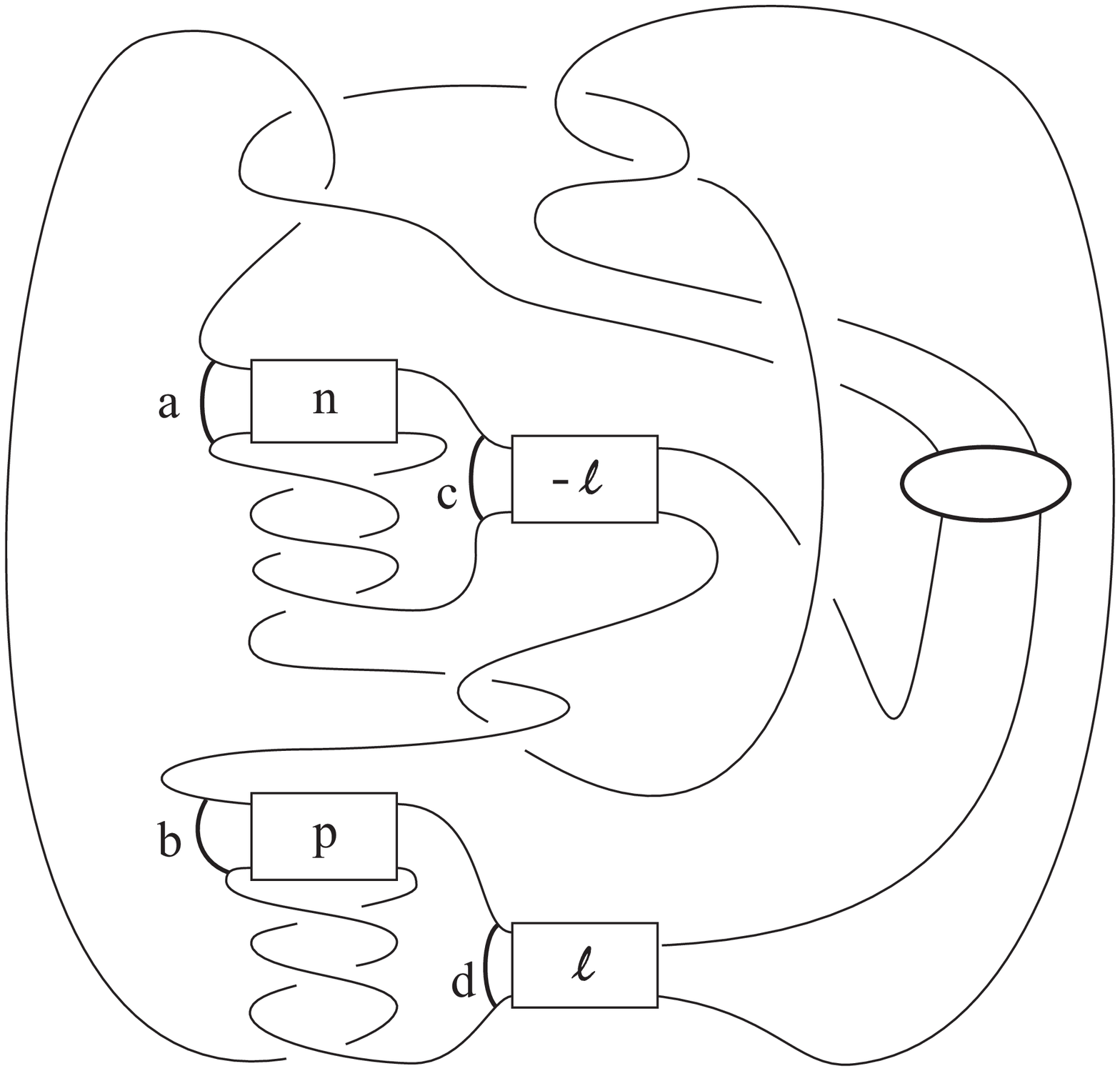}
\caption{Trivializable tangle $T(l, n, p)$: $n$ or $p$ is zero.}
\label{Tlnp}
\end{center}
\end{figure}

Let $T(l, n, p)$ be the tangle of Figure~\ref{Tlnp}, 
which is $T(A, B)$ in \cite[Figure~3(a)]{EM2} with
$A =R(n, -3, -l, 2, 1)$ and $B =R(p, -3, l)$.  
Then Lemma~4.1 in \cite{EM2} shows that 
$T(l, n, p)$ is a trivializable tangle
if and only if $n$ or $p$ is $0$. 
We summarize results in \cite[Section~4]{EM2} as follows. 

\begin{proposition}
\label{Tlnp+R}
\begin{enumerate}
\item The tangle $T(l, n, 0)$ has the following properties. 
	\begin{enumerate}
	\item
	$T(l, n, 0) + R(\infty)$ is a trivial knot. 

\smallskip	 
	 \item 
	$T(l, n, 0) + R(0)$ is the Montesinos--m link 
	$\displaystyle Mm(\frac{9l n - 3l + 1}{6l n - 2l - n + 1},\  l)$.
\smallskip	
	\item
	$T(l, n, 0) + R(1)$ is the Montesinos link 
	$\displaystyle M(3,\  \frac{9l n - 3l - 3n + 2}{-6l n + 2l + n - 1},\  
	\frac{l+1}{-l})$. 
	\end{enumerate}

\item The tangle $T(l, 0, p)$ has the following properties. 
\begin{enumerate}
	\item
	$T(l, 0, p) + R(\infty)$ is a trivial knot. 
	
	\item 
	$T(l, 0, p) + R(0)$ is the Montesinos--m link 
	$\displaystyle Mm( \frac{3l-1}{2l-1},\  \frac{3l p  - l - p}{3p-1})$. 
\smallskip
	\item
	$T(l, 0, p) + R(1)$ is the Montesinos link 
	$\displaystyle M(3,\  \frac{3l-2}{-2l+1},\  \frac{3l p - l + 2p - 1}{-3l p + l + p})$. 
	\end{enumerate}
\end{enumerate}
\end{proposition}

Let $\pi: S^3 \to S^3$ be the two--fold branched cover
along the trivial knot $T(l, n, p) +R(\infty)$,
where $n$ or $p$ is $0$.
Let $K(l, n, p)$ be the covering knot of the trivializable tangle $T(l, n, p)$, 
and $\gamma_{l, n, p}$ the covering slope
corresponding to $0$--untangle surgery on $T(l, n, p) + R(\infty)$, where $n$ or $p$ is 0.
Note that $1$--untangle surgery corresponds to 
$(\gamma_{l, n, p} - 1)$--surgery, where $\gamma_{l, n, p}
\in\mathbb{Q}$. 
We denote by $\mathcal{EM}\mathrm{I}$ 
the set of the Seifert surgeries $(K(l, n, p), \gamma_{l, n, p})$
and $(K(l, n, p), \gamma_{l, n, p}-1)$,
where $n$ or $p$ is 0. 
For brevity,
we often write $(K(l, n, p), \gamma)$ and
$(K(l, n, p), \gamma -1)$
for $(K(l, n, p), \gamma_{l, n, p})$
and $(K(l, n, p), \gamma_{l, n, p}-1)$,
respectively.

\begin{proposition}
[{\cite[Proposition~4.5]{EM2}}]    
\label{Klnp}
The knots $K(l,n,p)$ have the following Seifert surgeries. 

\begin{enumerate}

\item $\displaystyle K(l, n, 0)(\gamma_{l, n, 0})
= \mathbb{R}P^2(\frac{-6l n + 2l + n - 1}{9l n - 3l + 1},\ \frac{-1}{l})$,

\smallskip
\item $\displaystyle K(l,n,0)(\gamma_{l, n, 0}-1)
=S^2(\frac{-1}{3},\  \frac{6l n - 2l - n + 1}{9l n - 3l - 3n + 2},\ \frac{l}{l + 1})$,

\smallskip
\item $\displaystyle K(l,0,p)(\gamma_{l, 0, p})
=\mathbb{R}P^2(\frac{-2l + 1}{3l - 1},\  \frac{-3p + 1}{3l p  - l - p})$, and

\smallskip
\item $\displaystyle K(l,0,p)(\gamma_{l, 0, p}-1)
=S^2(\frac{-1}{3},\  \frac{2l - 1}{3l - 2},\  \frac{3l p - l - p}{3l p - l + 2p - 1})$.

\end{enumerate}
Furthermore, 
$\gamma_{l, n, 0} = 12l^2-4l-36l^2 n$ and 
$\gamma_{l, 0, p} =  12l^2-4l-4p(3l-1)$.
\end{proposition}

We detect seiferters for Seifert surgeries
given in Proposition~\ref{Klnp}
by applying Theorem~\ref{seiferter for covering knots}.

\begin{proposition}
\label{seiferter for Klnp}
\begin{enumerate}
\item
Let $c_a = \pi^{-1}(a)$,
where $a$ is the arc given by Figure~\ref{Tlnp} with $p = 0$. 
Then $c_a$ is a seiferter for both $(K(l, n, 0), \gamma)$ and 
$(K(l, n, 0), \gamma -1)$. 
\item
Let $c_b = \pi^{-1}(b)$,
where $b$ is the arc given by Figure~\ref{Tlnp} with $n = 0$. 
Then $c_b$ is a seiferter for both $(K(l, 0, p), \gamma)$ and 
$(K(l, 0, p), \gamma -1)$. 
\end{enumerate}
\end{proposition}

\begin{figure}[!ht]
\begin{center}
\includegraphics[width=0.85\linewidth]{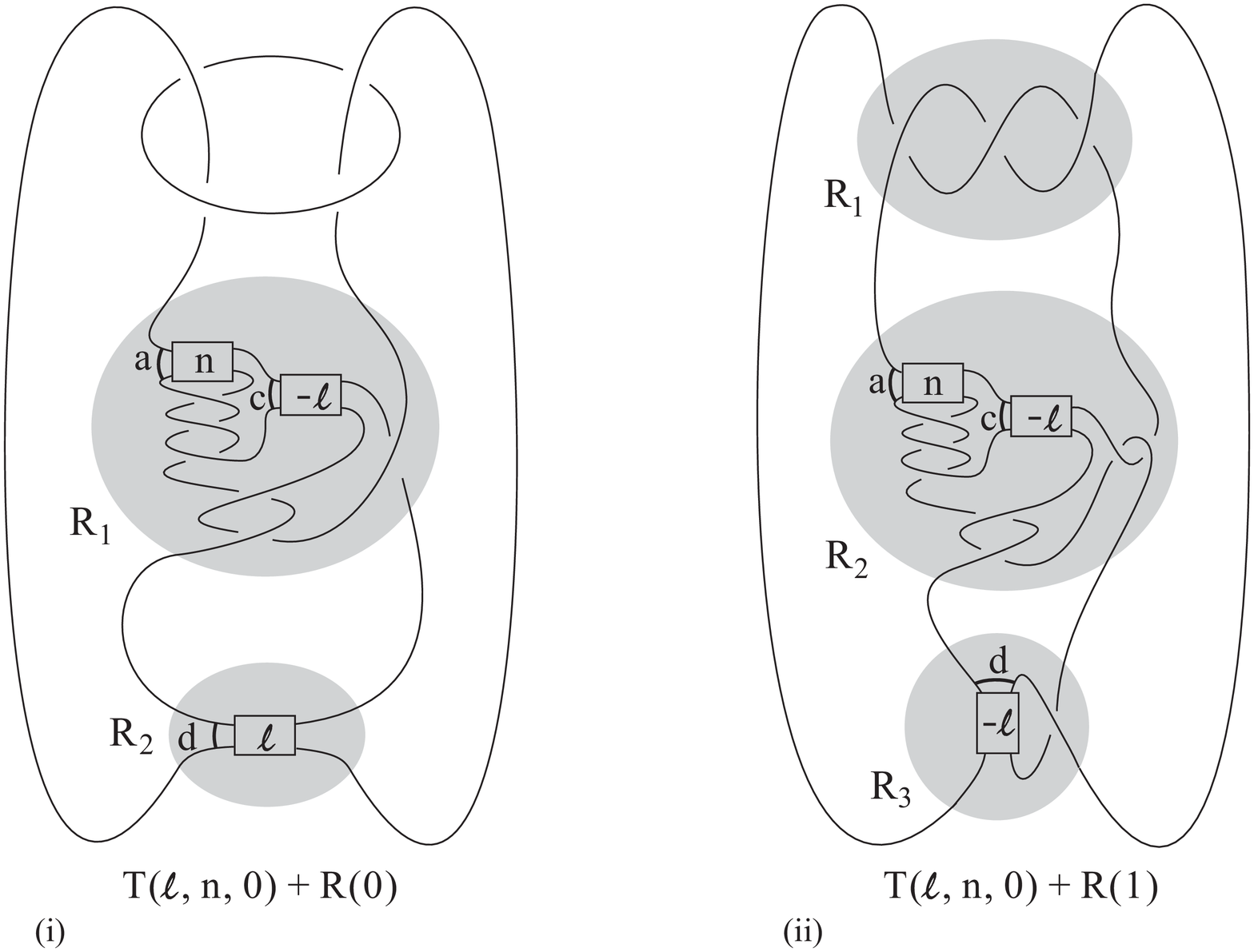}
\caption{}
\label{Tln0R01}
\end{center}
\end{figure}

\textsc{Proof of Proposition~\ref{seiferter for Klnp}.}
$(1)$ Isotoping the Montesinos--m link $T(l, n, 0) + R(0)$
together with the arcs $a, c, d$ in Figure~\ref{Tlnp},
we obtain Figure~\ref{Tln0R01}(i),
in which the Montesinos--m link is in a standard position.
In Figure~\ref{Tln0R01}(i),
the arc $a$ is a leading arc of the rational tangle 
$R_1 = R(\frac{9l n - 3l + 1}{6l n - 2l - n + 1})$.
Similarly,
isotoping the Montesinos link $T(l, n, 0) + R(1)$
together with the arcs $a, c, d$ in Figure~\ref{Tlnp},
we obtain Figure~\ref{Tln0R01}(ii).
In Figure~\ref{Tln0R01}(ii),
the Montesinos link is in a standard position,
and the arc $a$ is a leading arc of the rational tangle
$R_2 = R(\frac{9l n - 3l - 3n + 2}{-6l n + 2l + n - 1})$.

\begin{figure}[!ht]
\begin{center}
\includegraphics[width=0.35\linewidth]{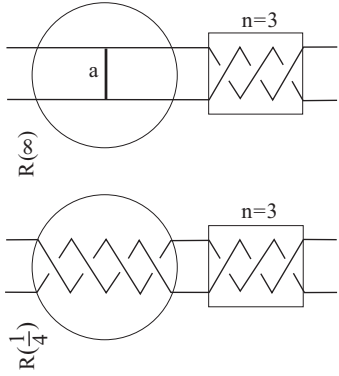}
\caption{$\displaystyle
\frac{1}{4}$--untangle surgery along $a$.}
\label{untangle}
\end{center}
\end{figure}

Apply $\frac{1}{n'}$--untangle surgery on
 $L = T(l, n, 0) + R(\infty)$
along the arc $a$.
More precisely, 
replace a trivial tangle $( N(a), N(a) \cap L )=R(\infty)$
with $R(\frac{1}{n'})$,
where $( N(a), N(a) \cap L )$ is identified with $R(\infty)$
by $(-\frac{\pi}{2})$--rotation about a line perpendicular to
the projection plane; see Figure~\ref{untangle} for example.  
We then obtain $T(l, n-n', 0) + R(\infty)$, 
which is a trivial knot by Proposition~\ref{Tlnp+R}(1)(i). 
Then, Theorem~\ref{seiferter for covering knots}(1) shows that 
$c_a$ is a seiferter for both $(K(l, n, 0), \gamma)$
and $(K(l, n, 0), \gamma -1)$. 
More precisely, 
by Lemma~\ref{fiber} $c_a$ is an exceptional fiber of index
$|9l n - 3l + 1|$ 
(resp.\ $|9l n - 3l - 3n + 2|$) 
in $K(l, n, 0)(\gamma)$
(resp.\ $K(l, n, 0)(\gamma -1)$).

\begin{figure}[!ht]
\begin{center}
\includegraphics[width=0.85\linewidth]{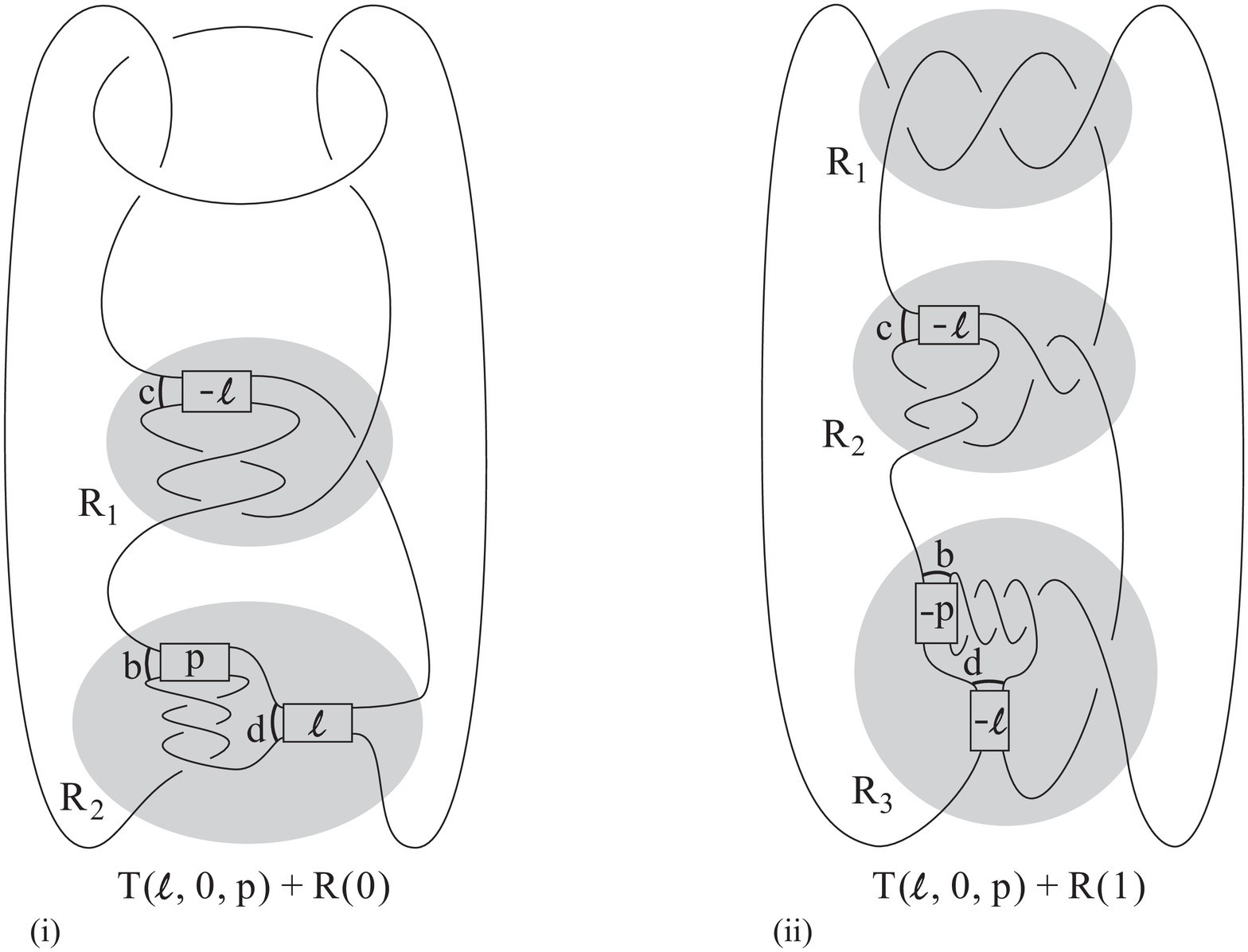}
\caption{}
\label{Tl0pR01}
\end{center}
\end{figure}

$(2)$ Let $b$ be the arc given in Figure~\ref{Tlnp}.  
Figures~\ref{Tl0pR01}(i) and (ii) give standard positions of 
$T(l, 0, p) + R(0)$ and $T(l, 0, p) + R(1)$,
respectively. 
Then $b$ is a leading arc of the rational tangle 
$R_2 = R(\frac{3l p - l - p}{3p - 1})$ 
(resp.\ $R_3
= R(\frac{3l p - l + 2p - 1}{-3l p + l + p})$) of 
$T(l, 0, p) + R(0)$ (resp.\ $T(l, 0, p) + R(1)$).  

Identify a trivial tangle $(N(b), N(b) \cap L)$,
where $L = T(l, 0, p) +R(\infty)$,
with $R(\infty)$ by $(-\frac{\pi}{2})$--rotation,
and perform $\frac{1}{p'}$--untangle surgery on $L$
along the arc $b$. 
Then we obtain $T(l, 0, p-p') + R(\infty)$,
which is a trivial knot by Proposition~\ref{Tlnp+R}(2)(i). 
It follows from Theorem~\ref{seiferter for covering knots}(1) that 
$c_b$ is a seiferter for both $(K(l, 0, p), \gamma)$ and $(K(l,0,p), \gamma -1)$. 
More precisely $c_b$ becomes an exceptional fiber of index $|3l p - l - p|$ 
(resp.\ $|3l p - l + 2p - 1|$) in $K(l, 0, p)(\gamma)$ 
(resp.\ $K(l,0,p)(\gamma -1)$). 
This establishes Proposition~\ref{seiferter for Klnp}. 
\hspace*{\fill} $\square$(Proposition~\ref{seiferter for Klnp})

\begin{remark}
\label{cd}
Let $c_d =\pi^{-1}( d )$ in $S^3$,
where $d$ is the arc as in Figure~\ref{Tlnp}. 
Since $d$ is a leading arc of the rational tangle
$R_2$ in the standard position of $T(l, n, 0) + R(0)$
in Figure~\ref{Tln0R01},
$c_d$ is a Seifert fiber in $K(l, n, 0)(\gamma)$.
However, $c_d$ is not a trivial knot in $S^3$,
so that $c_d$ is not a seiferter for $(K(l, n, 0), \gamma)$. 
Similarly, although $c_c = \pi^{-1}(c)$, where $c$
is the arc depicted in Figure~\ref{Tlnp}, is a fiber in 
$K(l, 0, p)(\gamma)$, it is not a trivial knot.
Hence $c_c$ is not a seiferter
for $(K(l, 0, p), \gamma)$. 
\end{remark}

Since $\frac{1}{n'}$--untangle surgery 
on $T(l, n, 0) + R(\infty)$
along the arc $a$ preserves the triviality
for any $n'$,
by Remark~\ref{framing} the untangle surgery
corresponds to $(-\frac{1}{n'})$--surgery (i.e.\ $n'$--twist) along the seiferter $c_a$. 
Similarly, 
$\frac{1}{p'}$--untangle surgery on $T(l, 0, p) + R(\infty)$
along the arc $b$
corresponds to $p'$--twist along the seiferter $c_b$. 
As observed in the proof of 
Proposition~\ref{seiferter for Klnp},
$\frac{1}{n'}$--untangle surgery on $T(l, n, 0)+R(\infty)$
along the arc $a$ yields $T(l, n-n', 0)+R(\infty)$. 
Since untangle surgeries along the arc $a$
do not affect the attached tangle $R(\infty)$
in $T(l, n, 0)+R(\infty)$,
the image of the covering slope $\gamma_{l, n, 0}$
under $n'$--twist along $c_a$
corresponds to replacing $R(\infty)$
in $T(l, n-n', 0) +R(\infty)$ with $R(0)$.
Thus, $n'$--twist along $c_a$ converts
$(K(l, n, 0), \gamma_{l, n, 0})$ to
$(K(l, n-n', 0), \gamma_{l, n-n', 0})$.
The same result holds for $p'$--twist along
the seiferter $c_b$.
Therefore, we obtain Proposition~\ref{Kl00} below.

\begin{proposition}
\label{Kl00}
\begin{enumerate}
\item
$n$--twist along the seiferter $c_a$ converts 
$(K(l, n, 0), \gamma_{l, n, 0})$ to 
$(K(l, 0, 0), \gamma_{l, 0, 0})$, 
and $(K(l, n, 0), \gamma_{l, n, 0} -1)$
to $(K(l, 0, 0), \gamma_{l, 0, 0} -1)$. 

\item
$p$--twist along the seiferter $c_b$ converts 
$(K(l, 0, p), \gamma_{l, 0, p})$ to 
$(K(l, 0, 0), \gamma_{l, 0, 0})$, 
and  $(K(l, 0, p), \gamma_{l, 0, p} -1)$
to $(K(l, 0, 0), \gamma_{l, 0, 0} -1)$. 
\end{enumerate}
\end{proposition}

Proposition~\ref{Kl00} gives paths from
the Seifert surgeries $(K(l, n, p), \gamma )$
(resp.\ $(K(l, n, p), \gamma -1)$)
to $(K(l, 0, 0), \gamma)$
(resp.\ $(K(l, 0, 0), \gamma)$),
where $n$ or $p$ is $0$.
Figure~\ref{lattice1} below shows
a portion of the Seifert Surgery Network; 
the horizontal lines (resp.\ the vertical lines) are generated by twisting along $c_a$ (resp.\ $c_b$).

\begin{figure}[!ht]
\begin{center}
\includegraphics[width=1.0\linewidth]{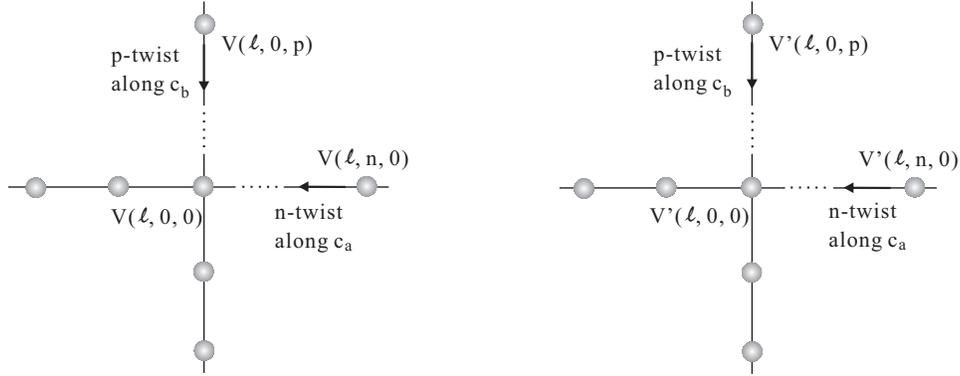}
\caption{
$V(l,n,p) = (K(l,n,p), \gamma_{l, n, p})$ and $V'(l,n,p) = (K(l,n,p), \gamma_{l, n, p}-1)$, 
where $n$ or $p$ is 0.
}
\label{lattice1}
\end{center}
\end{figure}

Let us show that there exist paths from 
$(K(l, 0, 0), \gamma)$ and
$(K(l, 0, 0), \gamma -1)$ 
to Seifert surgeries on torus knots. 

\begin{proposition}
\label{annular pair Kl00}
The pair $\{ c_c, c_d \}$ is an annular pair of seiferters for 
$(K(l, 0, 0), \gamma)$ and $(K(l, 0, 0), \gamma -1)$. 
\end{proposition} 

\begin{figure}[!ht]
\begin{center}
\includegraphics[width=0.85\linewidth]{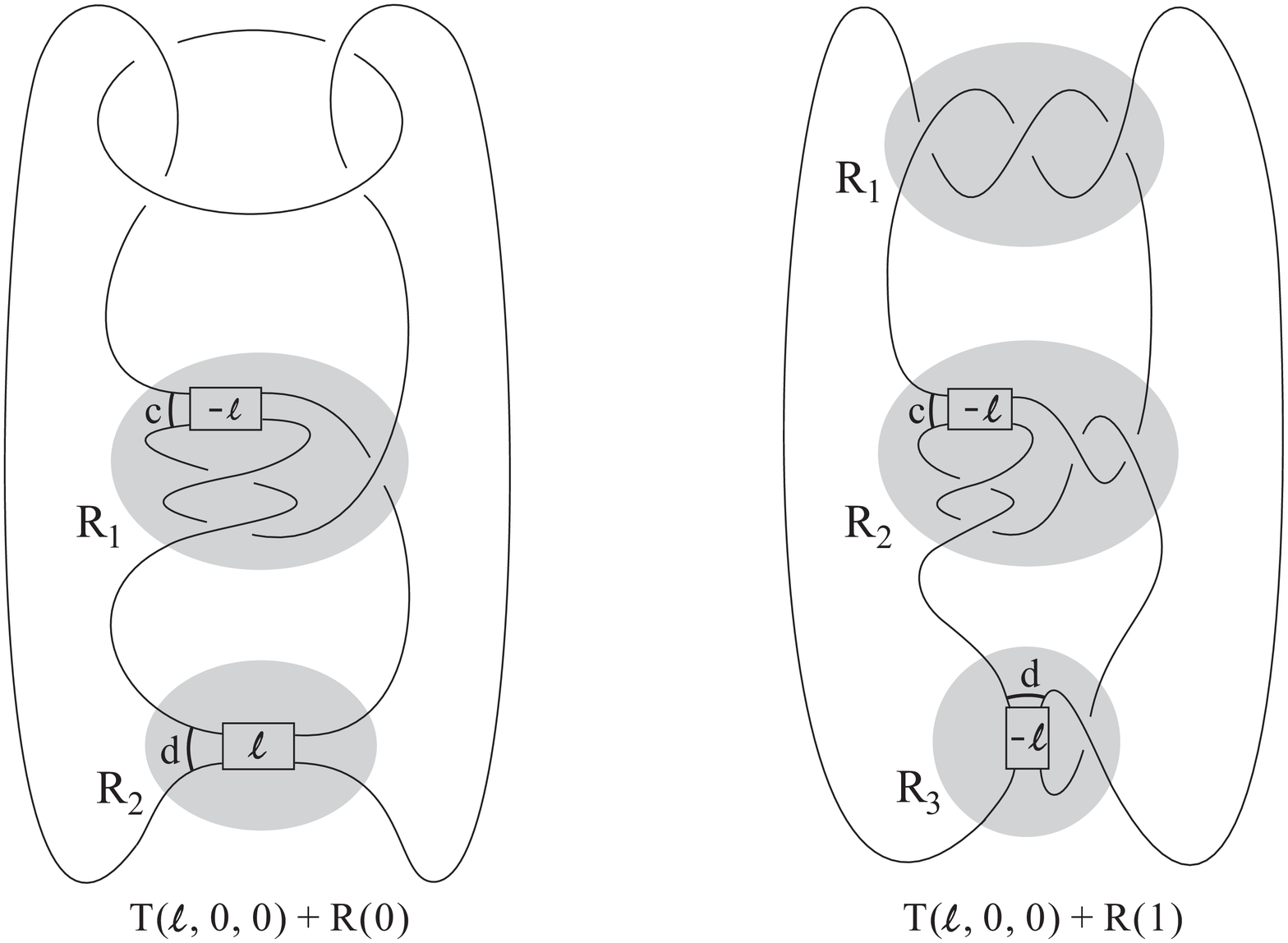}
\caption{}
\label{Tl00R01}
\end{center}
\end{figure}

\begin{figure}[!ht]
\begin{center}
\includegraphics[width=0.8\linewidth]{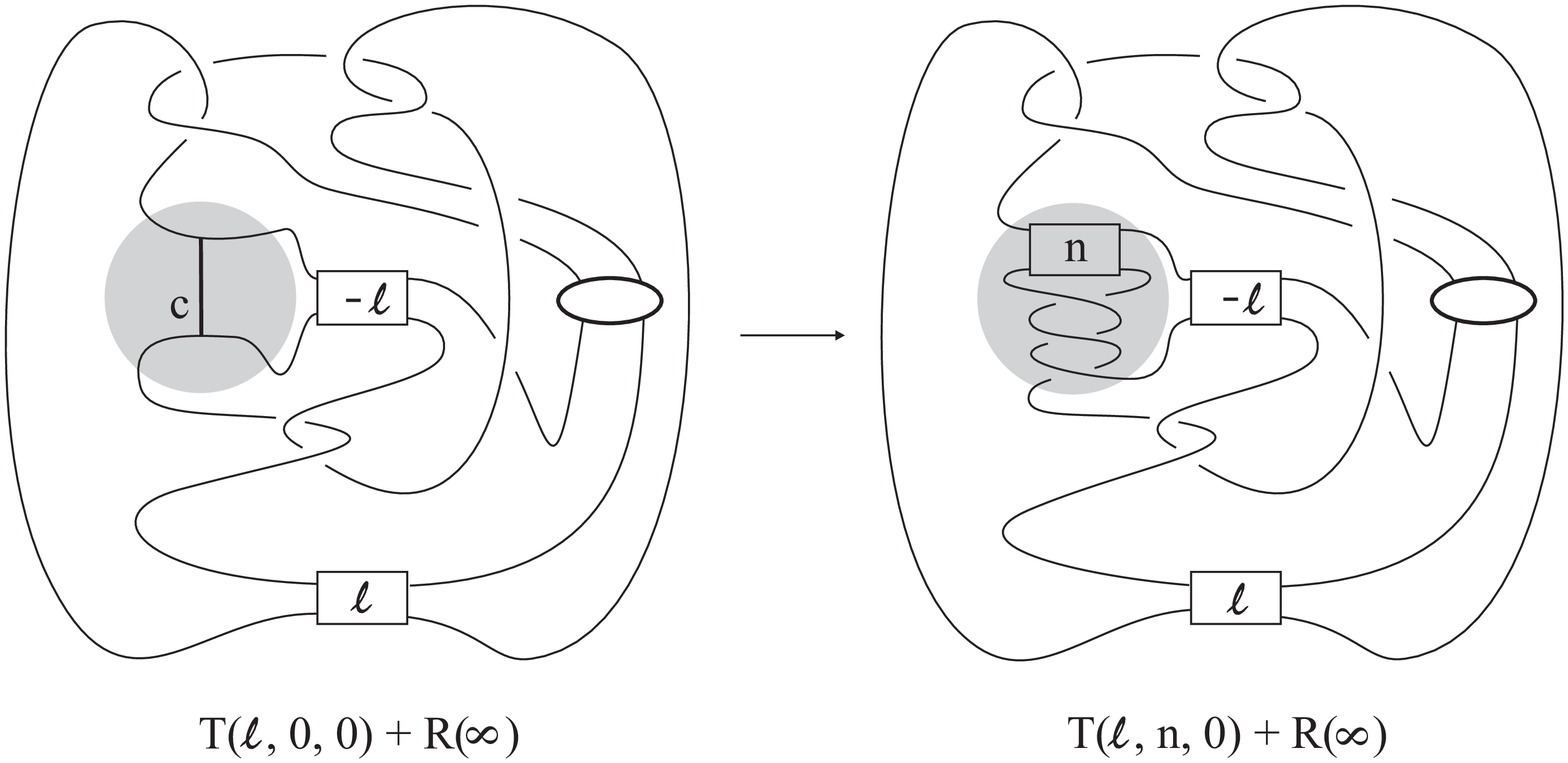}
\caption{}
\label{untangleTl00Tln0}
\end{center}
\end{figure}

\begin{figure}[!ht]
\begin{center}
\includegraphics[width=0.8\linewidth]{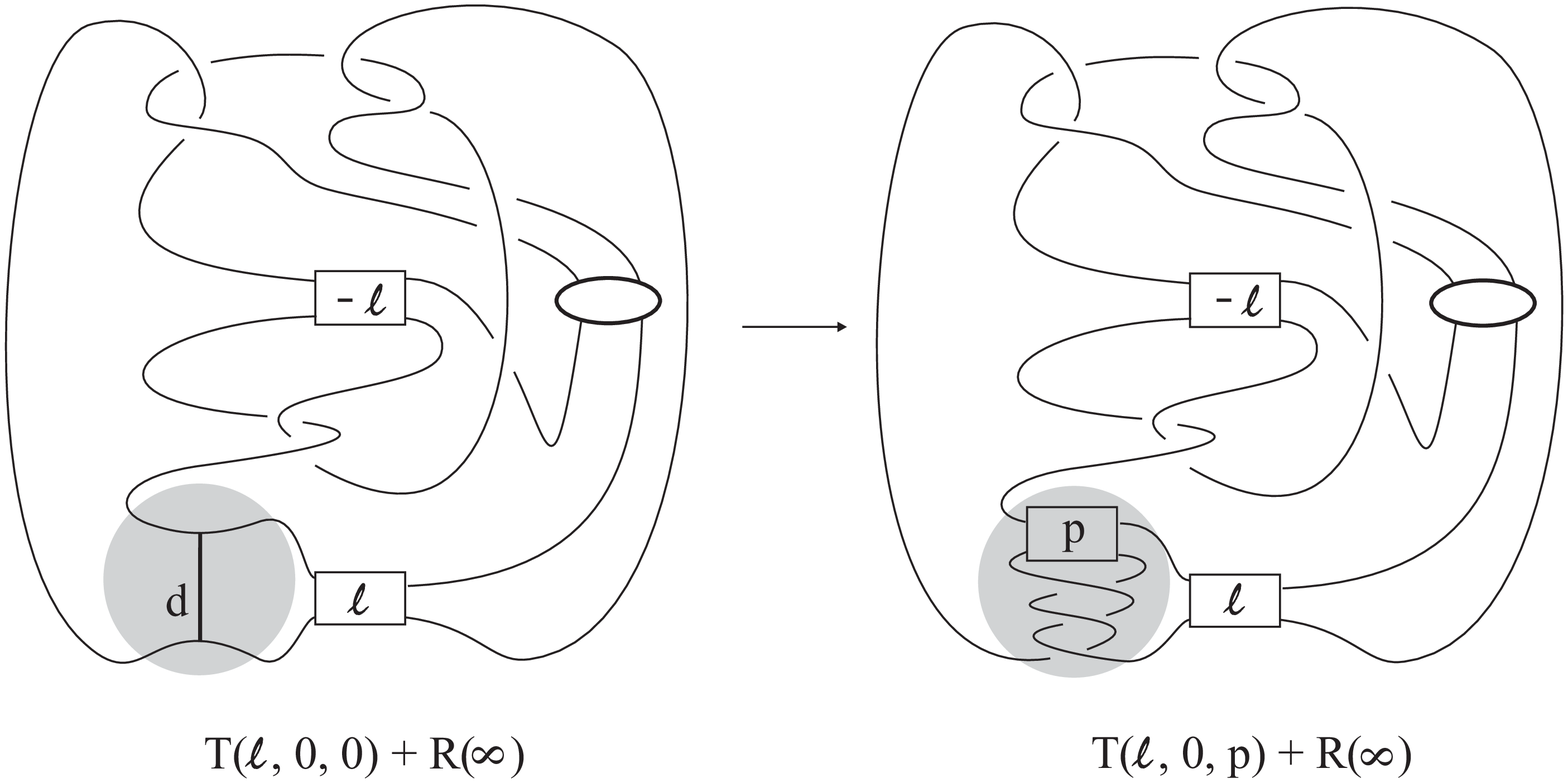}
\caption{}
\label{untangleTl00Tl0p}
\end{center}
\end{figure}

\textsc{Proof.}
We note that $c$ and $d$ are leading arcs of some $R_i$ in the standard positions of 
$T(l, 0, 0) + R(0)$ and $T(l, 0, 0) + R(1)$
given in Figure~\ref{Tl00R01}.  
Figure~\ref{untangleTl00Tln0} 
(resp. Figure~\ref{untangleTl00Tl0p}) 
shows that an untangle surgery along the arc $c$ (resp.\ $d$) 
converts the trivial knot $T(l, 0, 0) + R(\infty)$ 
to the trivial knot $T(l, n, 0) + R(\infty)$
(resp.\ $T(l, 0, p) + R(\infty)$). 
Theorem~\ref{annular pair for covering knots} then shows that
$\{c_c,\ c_d\}$ is a pair of seiferters for
$(K(l, 0, 0), \gamma)$
and $(K(l, 0, 0), \gamma -1)$. 

It remains to show that $c_c$ and $c_d$ cobound an annulus. 
Let $D$ be a rectangle as shown in Figure~\ref{Kl00D} below; 
Figure~\ref{Kl00D}(i) gives 
the part of $D$ which is on or above the projection plane,
(i.e.\ the upper part of $D$),
and (ii) gives the part of $D$ which is on or below the projection plane (i.e.\ the lower part of $D$).
The arcs $c$ and $d$ are opposite sides of $D$ and
the other sides of $D$ 
are contained in the trivial knot $T(l, 0, 0) + R(\infty)$.
Note that $D$ satisfies the condition of 
Theorem~\ref{annular pair for covering knots}. 
Hence, $c_c$ and $c_d$ cobound 
the annulus $\pi^{-1}(D)$. 
\hspace*{\fill} $\square$(Proposition~\ref{annular pair Kl00})

\begin{figure}[htbp]
\begin{center}
\includegraphics[width=1.0\linewidth]{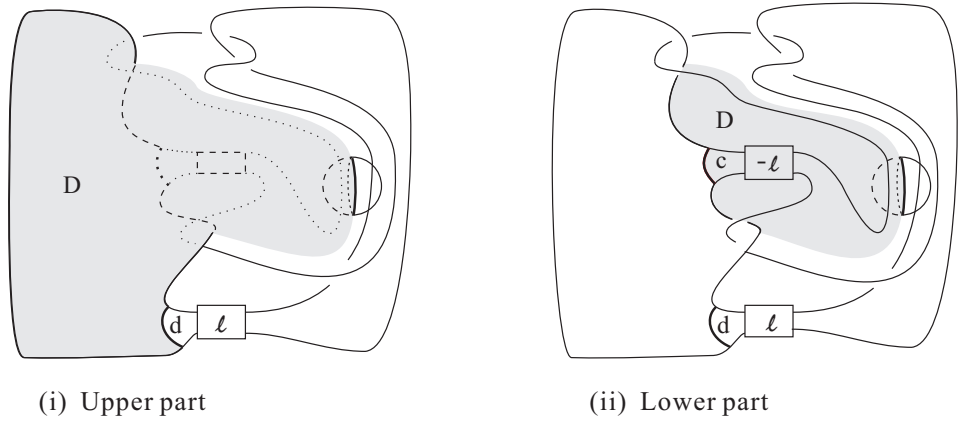}
\caption{A rectangle $D$ and  $T(l, 0, 0) + R(\infty)$.}
\label{Kl00D}
\end{center}
\end{figure}

\begin{lemma}
\label{K000}
$(-l)$-twist along the annular pair of seiferters 
$(c_c, c_d)$ converts $K(l, 0, 0)$ to
the trivial knot $K(0, 0, 0)$,
and $\gamma_{l,0,0}$ to $\gamma_{0,0,0} =0$.
\end{lemma}

\textsc{Proof.}
We see from Figure~\ref{Kl00D}(ii) that
$(\frac{-1}{l}, \frac{1}{l})$--untangle surgery 
along $(c, d)$ changes $L = T(l, 0, 0) + R(\infty)$ to 
$T(0, 0, 0) + R(\infty)$,
where $(N(c), N(c)\cap L)$ and $(N(d), N(d)\cap L)$
are identified with $R(\infty)$ as in Figure~\ref{untangle}.
Since the rectangle $D$ in Figure~\ref{Kl00D} intersects
the boundaries of $N(c)$ and $N(d)$ in their latitudes,
by Theorem~\ref{annular pair for covering knots}(2) 
the $(\frac{-1}{l}, \frac{1}{l})$--untangle surgery corresponds to 
$(-l)$--twist along $(c_c, c_d)$.
The triviality of $K(0, 0, 0)$ follows immediately by
observing that the union of $T(0, 0, 0)+R(\infty)$ 
and the spanning arc of $R(\infty)$ 
forms a $\theta$--curve standardly embedded in $S^3$ up to isotopy.   
Since untangle surgeries along $(c, d)$
do not affect the attached tangle
$R(\infty)$ in $T(l, 0, 0)+R(\infty)$,
the image of the surgery slope $\gamma_{l,0,0}$
corresponds to replacing $R(\infty)$ in $T(0, 0, 0)+R(\infty)$
with $R(0)$.
Hence, the image of $\gamma_{l,0,0}$ 
is $\gamma_{0,0,0}$, which is $0$ by
Proposition~\ref{Klnp}. 
\hspace*{\fill} $\square$(Lemma~\ref{K000})

\begin{figure}[htbp]
\begin{center}
\includegraphics[width=1.0\linewidth]{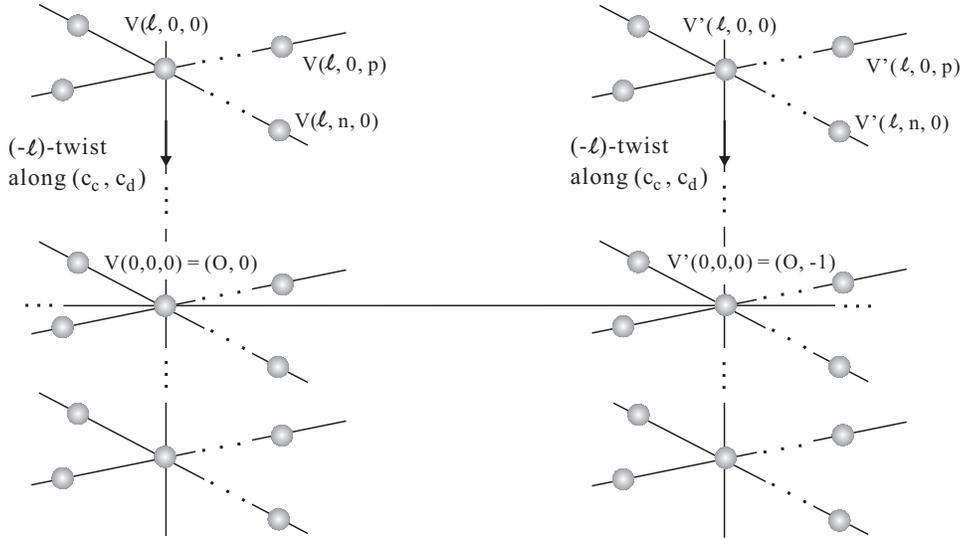}
\caption{
$\mathcal{EM}\mathrm{I}$ in the Seifert Surgery Network; 
$V(l,n,p) = (K(l,n,p), \gamma_{l, n, p})$ and $V'(l,n,p) = (K(l,n,p), \gamma_{l, n, p}-1)$, 
where $n$ or $p$ is 0.
}
\label{annularTree}
\end{center}
\end{figure}

\medskip
We now obtain Figure~\ref{annularTree} from
Proposition~\ref{Kl00} 
and Lemma~\ref{K000}.  
Each vertical line in Figure~\ref{annularTree} is generated by
twisting along $(c_c, c_d)$.
(The horizontal line in Figure~\ref{annularTree}
is generated by twisting along a meridian $\mu$ of the trivial knot $K(0,0,0)=O$;
$\mu$ is a seiferter for $(O, m)$ for any $m$.)
Figure~\ref{annularTree} gives explicit paths from 
$(K(l, n, p), \gamma)$ 
to $(K(0, 0, 0), 0) = (O, 0)$,
and from 
$(K(l, n, p), \gamma-1)$ to $(K(0, 0, 0), -1) = (O, -1)$,
where $n$ or $p$ is 0. 
\hspace*{\fill} $\square$(Theorem~\ref{connected}
for $\mathcal{EM}\mathrm{I}$) 

\medskip
In \cite[Example 9.25]{DMM1}
we discussed annular pairs of seiferters 
connecting $(K(l, 0, 0), \gamma_{l, 0, 0})$ and $(O, 0)$ from a different viewpoint. 

\subsection{The second family of Seifert surgeries
$\mathcal{EM}\mathrm{II}$}
\label{subsection:Klmnp}

\begin{figure}[!ht]
\begin{center}
\includegraphics[width=0.55\linewidth]{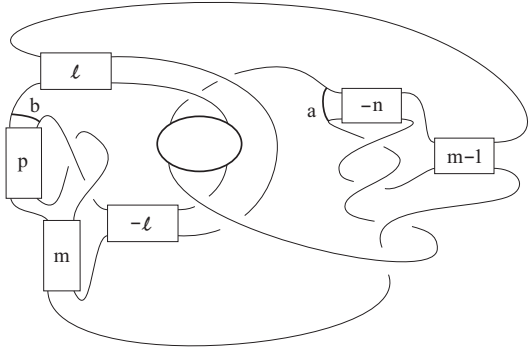}
\caption{Trivializable tangle $T(l, m, n, p)$, 
where $n = 0$ or $p = 0$.}
\label{Tlmnp}
\end{center}
\end{figure}

The tangle $T(l, m, n, p)$ of Figure~\ref{Tlmnp} 
is $\mathcal{B}(A, B, C)$ in Figure~9(a) of \cite[Section 5]{EM2} 
with $A = R(l)$, $B=R(p, -2, m, -l)$, $C = R(-n, 2, m-1, 2, 0)$; 
however, four crossings in Figure~9(a) should be revered.
Figure~\ref{Tlmnp} is the corrected diagram.
The second family $\mathcal{EM}\mathrm{II}$ consists of Seifert surgeries on 
the covering knots of $T(l, m, n, p)$.

\begin{proposition}
\label{Tlmnp+R}
\begin{enumerate}
\item
The tangle $T(l, m, n, 0)$ enjoys the following properties. 
	\begin{enumerate}
	\item
	$T(l, m, n, 0) + R(\infty)$ is a trivial knot.
	
\smallskip
	\item $T(l, m, n, 0) + R(0)$ is the Montesinos link 
	$\displaystyle M({l - 1},\  \frac{2mn - m - n + 1}{4mn - 2m + 1},\  \frac{l m + m - 1}{-m})$. 
	
\smallskip
	\item
	$T(l, m, n, 0) + R(1)$ is the Montesinos link 
	$\displaystyle M({l + 1},\  \frac{-2mn + m - n}{4mn - 2m + 1},\  \frac{- l m + m + 1}{m})$. 
	\end{enumerate}
	
\smallskip

\item
The tangle $T(l, m, 0, p)$ enjoys the following properties. 
	\begin{enumerate}
	\item
	$T(l, m, 0, p) + R(\infty)$ is a trivial knot.
	
\smallskip
	\item 
	$T(l, m, 0, p) + R(0)$ is the Montesinos link 
	\smallskip
	
	$\displaystyle M({l - 1},\ 
	\frac{2l mp - l m - l p + 2mp - m - 3p + 1}{-2mp + m + p},\  \frac{m - 1}{2m - 1})$. 
	
\smallskip
	\item
	$T(l, m, 0, p) + R(1)$ is the Montesinos link 
	\smallskip
	
	$\displaystyle M({l + 1},\   
	\frac{2l mp - l m - l p - 2mp + m - p + 1}{-2mp + m + p},\   \frac{-m}{2m-1})$. 
	\end{enumerate}
\end{enumerate}
\end{proposition}

\begin{figure}[htbp]
\begin{center}
\includegraphics[width=0.85\linewidth]{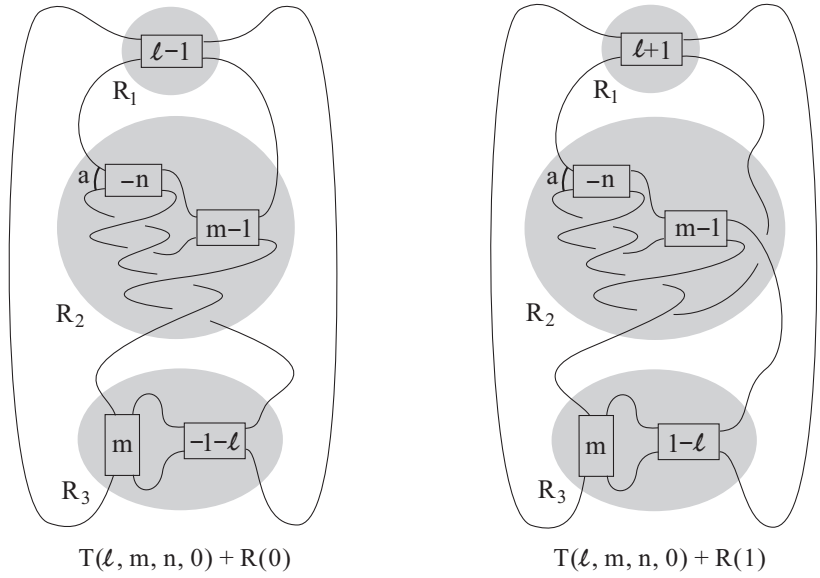}
\caption{Montesinos links $T(l, m, n, 0) + R(0)$,  
$T(l, m, n, 0) + R(1)$ in standard forms}
\label{Tlmn0R01}
\end{center}
\end{figure}

\begin{figure}[htbp]
\begin{center}
\includegraphics[width=0.85\linewidth]{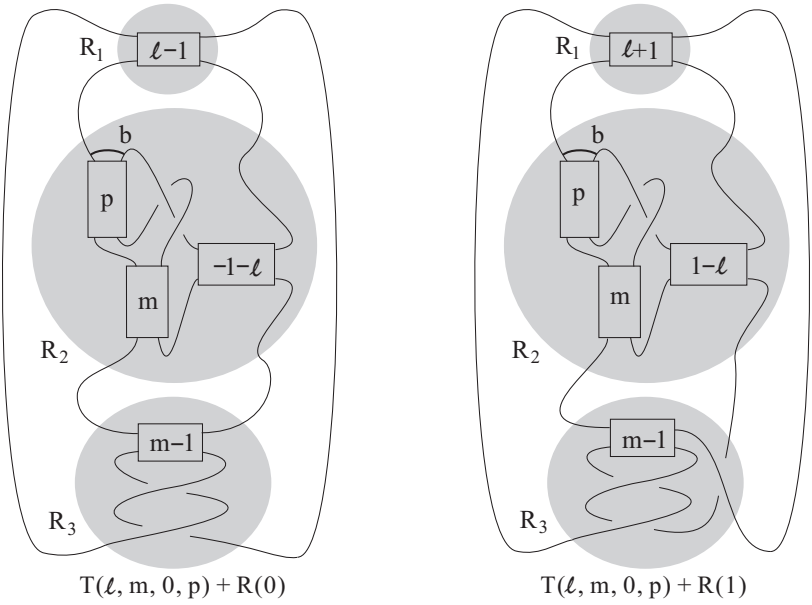}
\caption{Montesinos links $T(l, m, 0, p) + R(0)$,  
$T(l, m, 0, p) + R(1)$ in standard forms}
\label{Tlm0pR01}
\end{center}
\end{figure}

\textsc{Proof of Proposition~\ref{Tlmnp+R}.}
Assertions~(1)(i) and (2)(i) are straightforward. 
In fact, \cite[Lemma~5.1]{EM2} implies that $T(l, m, n, p)$
is a trivializable tangle if and only if $n$ or $p$ is 0. 
Assertions~(1)(ii) and (1)(iii) follow from
the standard positions of
$T(l, m, n, 0) + R(0)$ and $T(l, m, n, 0) + R(1)$ given in Figure~\ref{Tlmn0R01}. 
Similarly, $T(l, m, 0, p) + R(0)$ and $T(l, m, 0, p) + R(1)$ have the standard positions as in Figure~\ref{Tlm0pR01}.  
We then obtain (2)(ii) and (2)(iii). 
\hspace*{\fill} $\square$(Proposition~\ref{Tlmnp+R})

\medskip
Let $\pi: S^3 \to S^3$ be the two--fold cover branched along
$T(l, m, n, p) +R( \infty )$, where $n$ or $p$ is $0$.
Let $K(l, m, n, p)$ be the covering knot 
of the trivializable tangle $T(l, m, n, p)$,
and $\gamma_{l,m,n,p}$ the covering slope corresponding to 
$0$--untangle surgery on $T(l, m, n, p) + R(\infty)$,
where $n$ or $p$ is $0$. 
Then $1$--untangle surgery 
on $T(l, m, n, p) + R(\infty)$ 
corresponds to 
$(\gamma_{l, m, n, p} - 1)$--surgery on $K(l, m, n, p)$,
where $n$ or $p$ is $0$. 
We denote by $\mathcal{EM}\mathrm{II}$ 
the set of the Seifert surgeries 
$(K(l, m, n, p), \gamma_{l, m, n, p})$ and
$(K(l, m, n, p), \gamma_{l, m, n, p}-1)$,
where $n$ or $p$ is $0$.
For brevity,
we often write $(K(l,m, n, p), \gamma)$ and
$(K(l, m, n, p), \gamma -1)$
for $(K(l,m,n,p), \gamma_{l,m,n,p})$
and $(K(l,m,n,p), \gamma_{l,m,n,p}-1)$,
respectively.

Proposition~\ref{Tlmnp+R} can be translated into
the following assertion, 
which is a revision of Proposition~5.4 in \cite{EM2}. 

\begin{proposition}
\label{Klmnp}
$K(l,m,n,0)$ and $K(l,m,0,p)$ have the following Seifert surgeries. 
\begin{enumerate}

\item $\displaystyle K(l, m, n, 0)(\gamma_{l, m, n, 0}) =
S^2(\frac{-1}{l - 1},\  \frac{-4mn + 2m - 1}{2mn - m - n + 1},\  \frac{m}{l m + m - 1})$,

\smallskip
\item $\displaystyle K(l, m, n, 0)(\gamma_{l, m, n, 0}-1) =
S^2(\frac{-1}{l + 1},\  \frac{4mn - 2m + 1}{2mn - m + n},\  \frac{m}{l m - m - 1})$,

\smallskip
\item $\displaystyle K(l, m, 0, p)(\gamma_{l, m, 0, p}) =
S^2(\frac{-1}{l - 1},\  \frac{2mp - m - p}{2l mp - l m - l p + 2mp - m - 3p + 1},\  
\frac{-2m + 1}{m-1})$, and 

\smallskip
\item $\displaystyle K(l, m, 0, p)(\gamma_{l, m, 0, p}-1) =
S^2(\frac{-1}{l+1},\  \frac{2mp - m - p}{2l mp - 2mp -l m - l p + m - p + 1}, 
\frac{2m - 1}{m})$.
\end{enumerate}

Furthermore, 
$\gamma_{l, m, n, 0}=l(2m-1)(1-l m)+n(2l m-1)^2$, and
$\gamma_{l, m, 0, p}=l(2m-1)(1-l m)+p(2l m-l-1)^2$.
\end{proposition}

In order to detect seiferters 
for Seifert surgeries in Proposition~\ref{Klmnp}, 
as in the previous subsection, 
we apply Theorem~\ref{seiferter for covering knots}. 

\begin{proposition}
\label{seiferter for Klmnp}
\begin{enumerate}

\item
Let $c_a =\pi^{-1}(a)$, where $a$ is the arc given by Figure~\ref{Tlmnp} with $p = 0$. 
Then $c_a$ is a seiferter for both $(K(l, m, n, 0), \gamma)$ and 
$(K(l, m, n, 0), \gamma -1)$. 

\item
Let $c_b = \pi^{-1}(b)$, where $b$ is the arc given by Figure~\ref{Tlmnp} with $n = 0$. 
Then $c_b$ is a seiferter for both $(K(l, m, 0, p), \gamma)$ and 
$(K(l, m, 0, p), \gamma -1)$. 

\end{enumerate}
\end{proposition}

\textsc{Proof of Proposition~\ref{seiferter for Klmnp}.}
$(1)$ Figure~\ref{Tlmn0R01} gives standard positions of 
the following Montesinos links:

$${T(l, m, n, 0) + R(0) = M({l - 1},\  \frac{2mn - m - n + 1}{4mn - 2m + 1},\  \frac{l m + m - 1}{-m})}$$ 

and 

$${T(l, m, n, 0) + R(1) = M({l + 1},\  \frac{-2mn + m - n}{4mn - 2m + 1},\  \frac{-l m + m + 1}{m})}.$$
\noindent
Figure~\ref{Tlmn0R01} also shows that the arcs $a$ are leading arcs of $R_2$ in standard positions of these Montesinos links. 

Apply $\frac{1}{n'}$--untangle surgery on $T(l, m, n, 0) + R(\infty)$ along $a$ as in Figure~\ref{untangle}. 
We then obtain $T(l, m, n-n', 0) + R(\infty)$, 
which is a trivial knot by Proposition~\ref{Tlmnp+R}(1)(i). 
Now Theorem~\ref{seiferter for covering knots}(1) shows that 
$c_a$ is a seiferter for both $(K(l, m, n, 0), \gamma)$ and $(K(l, m, n, 0), \gamma-1)$.  
More precisely, 
by Lemma~\ref{fiber} $c_a$ is an exceptional fiber of index $|2mn - m - n + 1|$ 
(resp.\ $|2mn - m + n|$) 
in $K(l, m, n, 0)(\gamma)$ 
(resp.\ $K(l, m, n, 0)(\gamma -1)$).

$(2)$ Figure~\ref{Tlm0pR01} gives standard positions of 
the following Montesinos links: 

$${T(l, m, 0, p) + R(0) = M({l - 1},\  \frac{2l mp - l m - l p + 2mp - m - 3p + 1}{-2mp + m + p},\  
\frac{m - 1}{2m - 1})}$$ 

and 

$${T(l, m, 0, p) + R(1) = M({l + 1},\  \frac{2l mp  - l m - l p - 2mp + m - p + 1}{-2mp + m + p},\  
\frac{-m}{2m - 1})}.$$
\noindent
Figure~\ref{Tlm0pR01} also shows that the arcs $b$ are leading arcs of $R_2$ in standard positions of these Montesinos links.

Note that $\frac{1}{p'}$-untangle surgery on $T(l, m, 0, p) + R(\infty)$ 
along the arc $b$ as in Figure~\ref{untangle} yields the trivial knot 
$T(l, m, 0, p-p') + R(\infty)$ (Proposition~\ref{Tlmnp+R}(2)(i)). 
It follows from Theorem~\ref{seiferter for covering knots}(1) that 
$c_b$ is a seiferter for both $(K(l, m, 0, p), \gamma)$ and $(K(l, m, 0, p), \gamma-1)$; 
$c_b$ becomes an exceptional fiber of index 
$|2l mp - l m - l p + 2mp - m - 3p + 1|$ 
(resp.\  $|2l mp - l m - l p  -2mp + m - p + 1|$) in $K(l, m,0,p)(\gamma)$ 
(resp.\ $K(l, m, 0, p)(\gamma-1)$). 
This establishes Proposition~\ref{seiferter for Klmnp}. 
\hspace*{\fill} $\square$(Proposition~\ref{seiferter for Klmnp}) 

\medskip

In \cite{MM3} we demonstrated that $c_a$ is a seiferter for 
$(K(l, m, n, 0), \gamma)$
by applying a similar observation, 
but we checked the triviality of $c_a$
by drawing an diagram of $c_a$.

Since $\frac{1}{n'}$--untangle surgery on $T(l, m, n, 0) + R(\infty)$ along the arc $a$
preserves the triviality for any $n'$, 
by Remark~\ref{framing}
the untangle surgery
corresponds to $(-\frac{1}{n'})$--surgery (i.e.\ $n'$--twist) along the seiferter $c_a$.
Similarly, 
$\frac{1}{p'}$--untangle surgery on $T(l, m, 0, p) + R(\infty)$ along $b$
corresponds to $p'$--twist along the seiferter $c_b$. 
Since untangle surgeries along the arc $a$
do not affect the attached tangle $R(\infty)$
in $T(l, m, n, 0) + R(\infty)$,
the image of the covering slope $\gamma_{l, m, n, 0}$
under $n'$--twist along $c_a$
corresponds to replacing $R(\infty)$
in $T(l, m, n-n', 0) +R(\infty)$ with $R(0)$.
Thus, $n'$--twist along $c_a$ converts
$(K(l, m, n, 0), \gamma_{l, m, n, 0})$ to
$(K(l, m, n-n', 0), \gamma_{l, m, n-n', 0})$.
The same result holds for $p'$--twist along $c_b$.
Therefore, we obtain Proposition~\ref{Klm00}.

\begin{proposition}
\label{Klm00}
\begin{enumerate}
\item
$n$--twist along the seiferter $c_a$ converts 
$(K(l, m, n, 0), \gamma_{l, m, n, 0})$ to $(K(l, m, 0, 0), \gamma_{l, m, 0, 0})$,
and $(K(l, m, n, 0), \gamma_{l, m, n, 0} -1)$ to $(K(l, m, 0, 0),\gamma_{l, m, 0, 0}-1)$.  
\item
$p$--twist along the seiferter $c_b$ converts 
$(K(l, m, 0, p), \gamma_{l, m, 0, p})$ to 
$(K(l, m, 0, 0), \gamma_{l, m, 0, 0})$,
and $(K(l, m, 0, p), \gamma_{l, m, 0, p} -1)$ to $(K(l, m, 0, 0),\gamma_{l, m, 0, 0}-1)$.
\end{enumerate}
\end{proposition}

\begin{figure}[!ht]
\begin{center}
\includegraphics[width=1.0\linewidth]{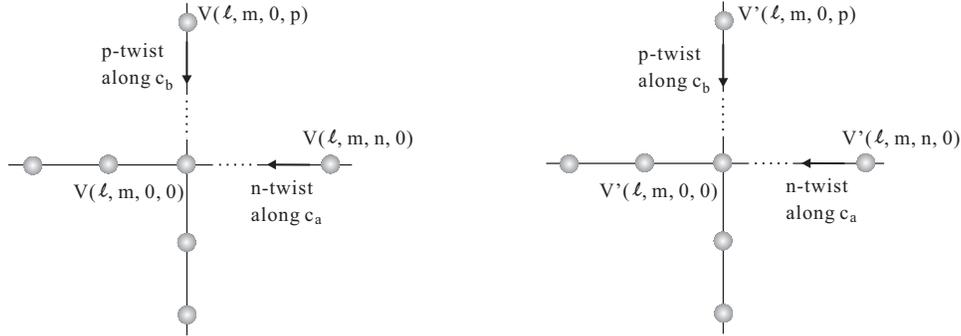}
\caption{$V(l, m, n, p) = (K(l, m, n, p), \gamma_{l,m,n,p})$ and 
$V'(l, m, n, p) = (K(l, m, n, p), \gamma_{l,m,n,p}-1)$, 
where $n$ or $p$ is 0.}
\label{lattice3}
\end{center}
\end{figure}

Proposition~\ref{Klm00}(1) and (2) give horizontal lines and 
vertical lines in Figure~\ref{lattice3}, respectively.  

\begin{proposition}
\label{isotopyT}
\begin{enumerate}
\item
The tangle $T(l, m, 0, 1)$ is equivalent to
$T(l, m-1, 1, 0)$.
\item
$(K(l, m, 0, 1), \gamma_{l, m, 0, 1}) 
= (K(l, m-1, 1, 0), \gamma_{l, m-1, 1, 0})$.  
\end{enumerate}
\end{proposition}

\begin{figure}[htbp]
\begin{center}
\includegraphics[width=1.0\linewidth]{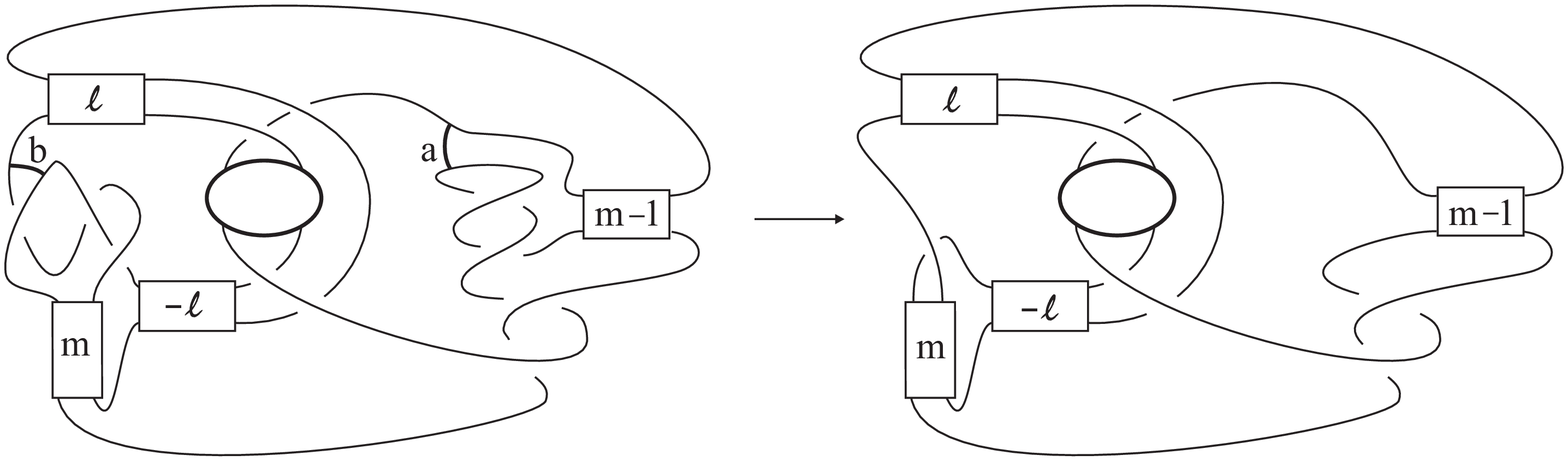}
\caption{An isotopy of $T(l, m, 0, 1)$.}
\label{Tlm01}
\end{center}
\end{figure}

\begin{figure}[htbp]
\begin{center}
\includegraphics[width=1.0\linewidth]{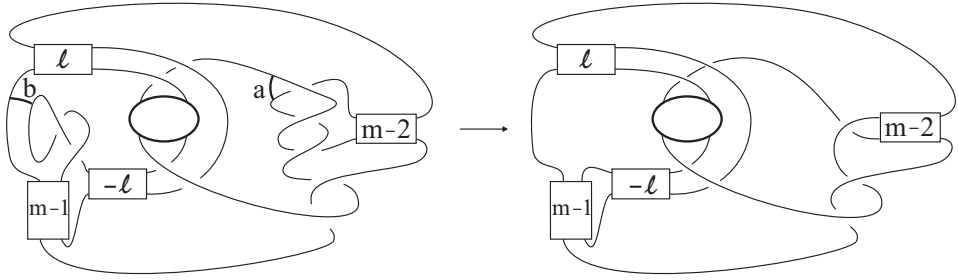}
\caption{An isotopy of $T(l, m-1, 1, 0)$.}
\label{Tlm-110}
\end{center}
\end{figure}

\textsc{Proof of Proposition~\ref{isotopyT}.}
$(1)$ We give a pictorial proof.  
The tangles on the right hand sides of
Figures~\ref{Tlm01} and \ref{Tlm-110} are equivalent.
Hence, the isotopies in Figures~\ref{Tlm01} and \ref{Tlm-110}
show that $T(l, m, 0, 1)$ and $T(l, m-1, 1, 0)$ are
equivalent to the same tangle, and thus equivalent to each other.

$(2)$ 
Assertion~(1) shows that
there is an automorphism of $S^3$ which sends
$T(l, m, 0, 1)+R(\infty)$ to $T(l, m-1, 1, 0)+R(\infty)$
and fixes $R(\infty)$.
Hence, there is an automorphism of $S^3$
which sends the covering knot $K(l,m ,0, 1)$
to $K(l, m-1, 1, 0)$, and the covering slope
$\gamma_{l, m, 0, 1}$ to $\gamma_{l, m-1, 1, 0}$.
This completes the proof. 
(We can obtain
$\gamma_{l, m, 0, 1} = \gamma_{l, m-1, 1, 0}$ 
directly from the formula of 
covering slopes in Proposition~\ref{Klmnp}.) 
\hspace*{\fill} $\square$(Proposition~\ref{isotopyT})

\medskip

\begin{figure}[htbp]
\begin{center}
\includegraphics[width=0.9\linewidth]{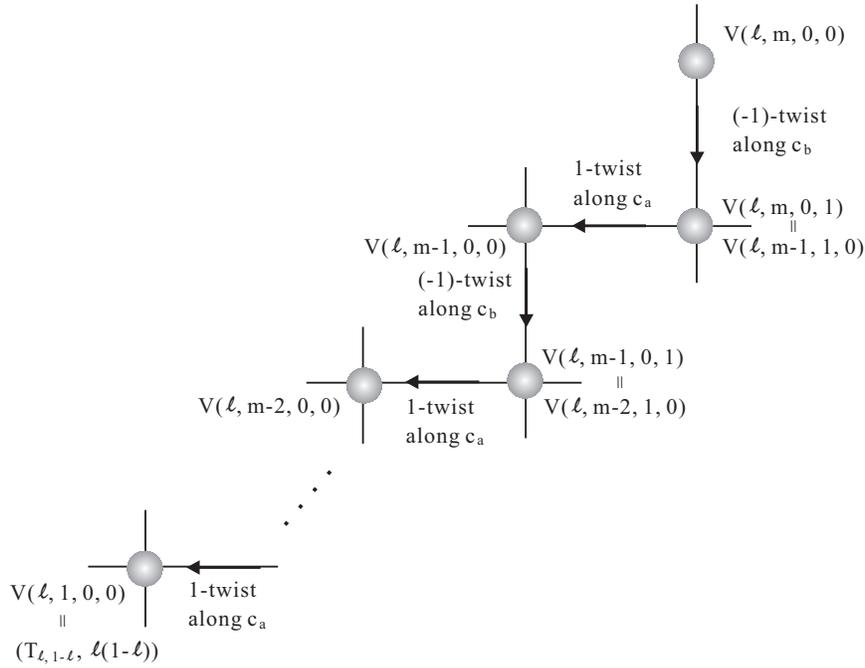}
\caption{$V(l, m, n, p) = (K(l, m, n, p), \gamma_{l,m,n,p})$}
\label{Klm00lattice}
\end{center}
\end{figure}

Applying Propositions~\ref{Klm00}
and \ref{isotopyT}(2) repeatedly, 
we find a path 
from $(K(l, m, 0, 0), \gamma)$ to 
$( K(l, 1, 0, 0), \gamma )$
as in Figure~\ref{Klm00lattice},
and a path
from $(K(l, m, 0, 0), \gamma-1)$ to 
$( K(l, 1, 0, 0), \gamma -1)$. 
Claim~\ref{Kl100=torus knot} below shows that
$K(l, 1, 0, 0)$ is a torus knot.
\hspace*{\fill} $\square$(Theorem~\ref{connected}
for $\mathcal{EM}\mathrm{II}$)

\begin{claim}
\label{Kl100=torus knot}
$(K(l, 1, 0, 0), \gamma_{l, 1, 0, 0}) = (T_{l, 1-l},\ l(1-l))$.  
\end{claim}

\begin{figure}[htbp]
\begin{center}
\includegraphics[width=0.65\linewidth]{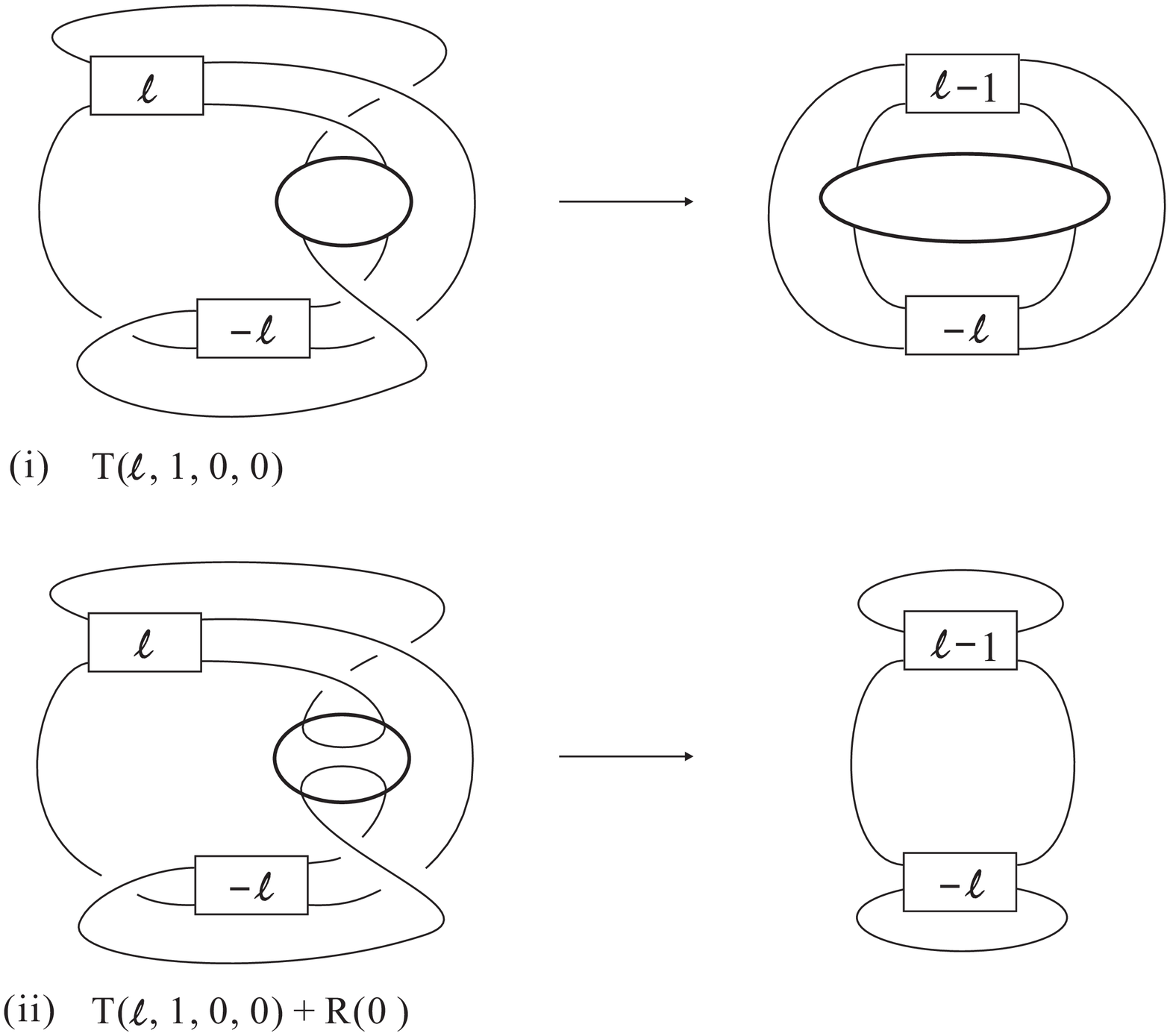}
\caption{}
\label{Tl100}
\end{center}
\end{figure}

\textsc{Proof.}
The tangle $(B, t)=T(l, 1, 0, 0)$ is a partial sum of
two rational tangles $R(l -1)$ and $R(-l)$;
see Figure~\ref{Tl100}.
Thus the exterior of $K(l, 1, 0, 0)$,
which is the two--fold cover of $B$ along $t$,
is a Seifert fiber space with two exceptional fibers of
indices $|l|, |l-1|$.
It follows that $K(l, 1, 0, 0)$ is
an $(l, 1-l)$ or $(l, l -1)$ torus knot. 
Figure~\ref{Tl100} shows that 
$T(l, 1, 0, 0) +R(0)$ is a connected sum of two $2$--bridge links, 
so that $K(l, 1, 0, 0)(\gamma_{l, 1, 0, 0})$ is a connected sum of 
two lens spaces.
Referring to Proposition~\ref{Klmnp}, 
we see that the the reducing slope $\gamma_{l, 1, 0, 0}$
equals $l(1-l)$, 
Thus $K(l, 1, 0, 0)$ is the $(l, 1-l)$ torus knot. 
\hspace*{\fill} $\square$(Claim~\ref{Kl100=torus knot})

\subsection{The third family of Seifert surgeries
$\mathcal{EM}\mathrm{III}$}
\label{subsection:KABCi}

The third family $\mathcal{EM}\mathrm{III}$ consists of
Seifert surgeries on knots $K(A, B, C)$,
which are the covering knots of the tangles $Q(A, B, C)$ below. 

\begin{figure}[htbp]
\begin{center}
\includegraphics[width=0.4\linewidth]{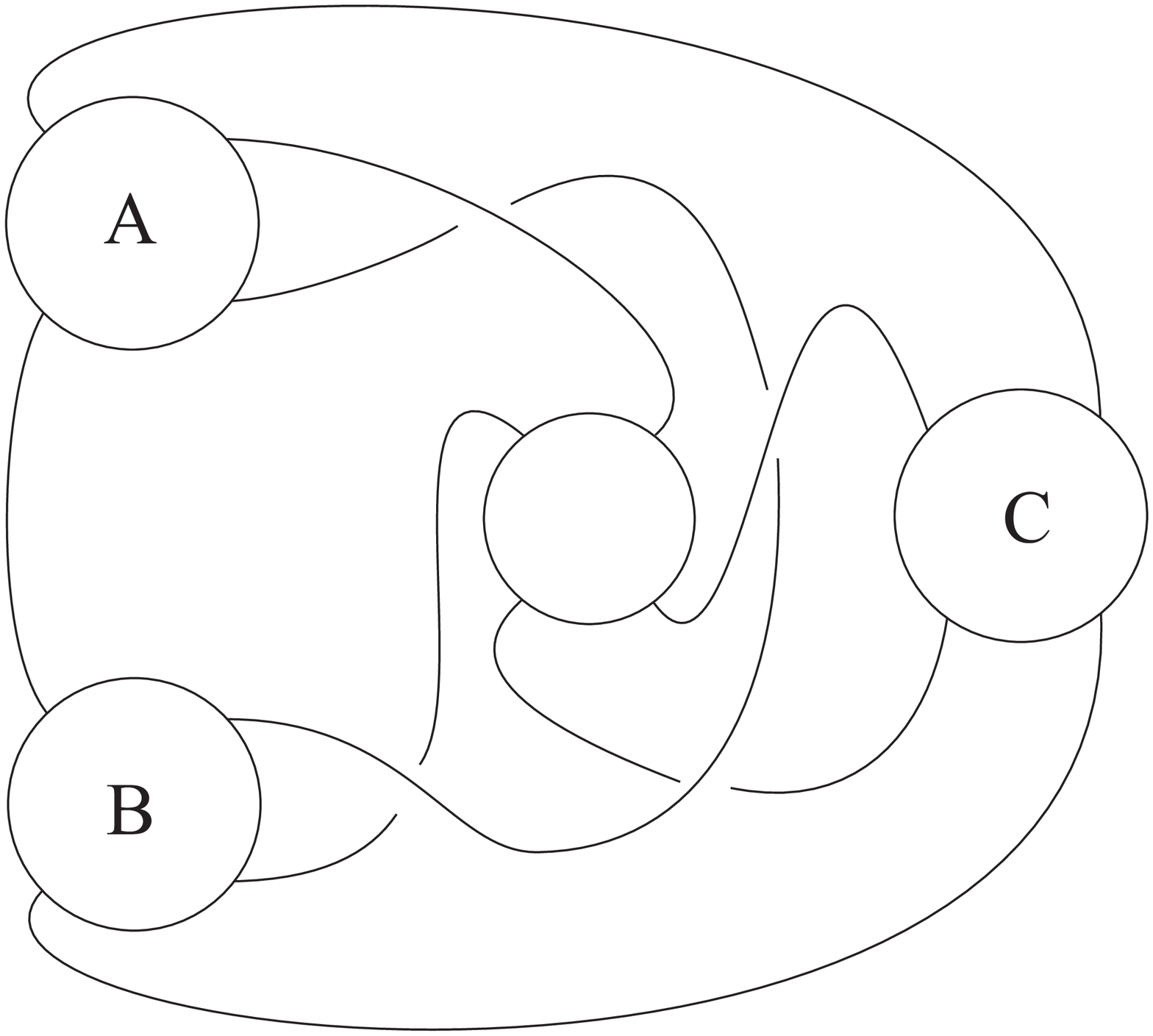}
\caption{Tangle $Q(A, B, C)$,
where $A, B, C$ are rational tangles.}
\label{Qabc}
\end{center}
\end{figure}

Let $Q(A, B, C)$ be the tangle of Figure~\ref{Qabc}, 
where $A, B, C$ are rational tangles;
it is also denoted by
$Q({\frac{\alpha_1}{\beta_1},  
\frac{\alpha_2}{\beta_2}}, {\frac{\alpha_3}{\beta_3}})$
if $A, B, C$ correspond to rational numbers
${\frac{\alpha_1}{\beta_1},  
\frac{\alpha_2}{\beta_2}}, {\frac{\alpha_3}{\beta_3}}$, respectively. 
This tangle was studied by the second author in \cite{EM2}
to produce an infinite family of 
Seifert surgeries on hyperbolic knots. 
Assume that 
$\frac{\alpha_i}{\beta_i} \ne \infty, 0, 1, 2$ for $i = 1, 2$ and 
$\frac{\alpha_3}{\beta_3} \ne \infty, 0, \pm 1, -\frac{1}{2}, 2$, 
for otherwise 
$Q(\frac{\alpha_1}{\beta_1}, \frac{\alpha_2}{\beta_2}, \frac{\alpha_3}{\beta_3})$ 
is a trivial tangle or a Montesinos tangle, 
i.e. a partial sum of rational tangles. 

\begin{proposition}
[{\cite[Lemmas 6.1 and 6.2]{EM2}}]  
\label{TABC}
$Q=Q(\frac{\alpha_1}{\beta_1}, \frac{\alpha_2}{\beta_2}, 
\frac{\alpha_3}{\beta_3})$ 
has the following properties up to interchanging $\frac{\alpha_1}{\beta_1}$ and 
$\frac{\alpha_2}{\beta_2}$.  

\begin{enumerate}
\item
$Q+ R(\infty)$ 
is a trivial knot if and only if 
$(\mathrm{i})$ or $(\mathrm{ii})$ below holds. 
	\begin{enumerate}
	\item
	For some integer $n$, $\frac{\alpha_3}{\beta_3} = \frac{1}{n}$ and $n\alpha_1\alpha_2 + \alpha_1\beta_2 + \beta_1\alpha_2 = \pm 1$. 
	\item
	For some integer $p$, $\frac{\alpha_1}{\beta_1} = \frac{1}{p}$ and $p\alpha_2\alpha_3 + \alpha_2\beta_3 + \beta_2\alpha_3 = \pm 1$. 
	\end{enumerate}

\item 
$Q+ R(-1)$ is the Montesinos link 
${M(\frac{\alpha_1 - 2\beta_1}{\beta_1}, 
\frac{\alpha_2 - 2\beta_2}{\beta_2}, \frac{\alpha_3+ \beta_3}{\beta_3})}$.
\end{enumerate}

\end{proposition}
 
Let us denote by 
$K(A, B, C)$ the covering knot of the trivializable tangle 
$Q(A, B, C)$, where $A,B,C$ satisfy condition~(i) or (ii) in Proposition~\ref{TABC}(1).  
Most of $K(A, B, C)$ are hyperbolic knots as shown in \cite{EM2}.  
Let $\gamma$ be the covering slope corresponding
to $(-1)$--untangle surgery on $Q(A, B, C) + R(\infty)$. 
Then Proposition~\ref{TABC}(2) shows that 
$K(\frac{\alpha_1}{\beta_1}, \frac{\alpha_2}{\beta_2}, \frac{\alpha_3}{\beta_3})(\gamma)$ is a Seifert fiber space over $S^2$ 
with three exceptional fibers of indices $|\alpha_1 - 2\beta_1|, |\alpha_2 - 2\beta_2|$, 
$|\alpha_3 + \beta_3|$. 
We denote by $\mathcal{EM}\mathrm{III}$ 
the set of the Seifert surgeries $(K(A, B, C), \gamma)$.

\begin{figure}[htbp]
\begin{center}
\includegraphics[width=1.0\linewidth]{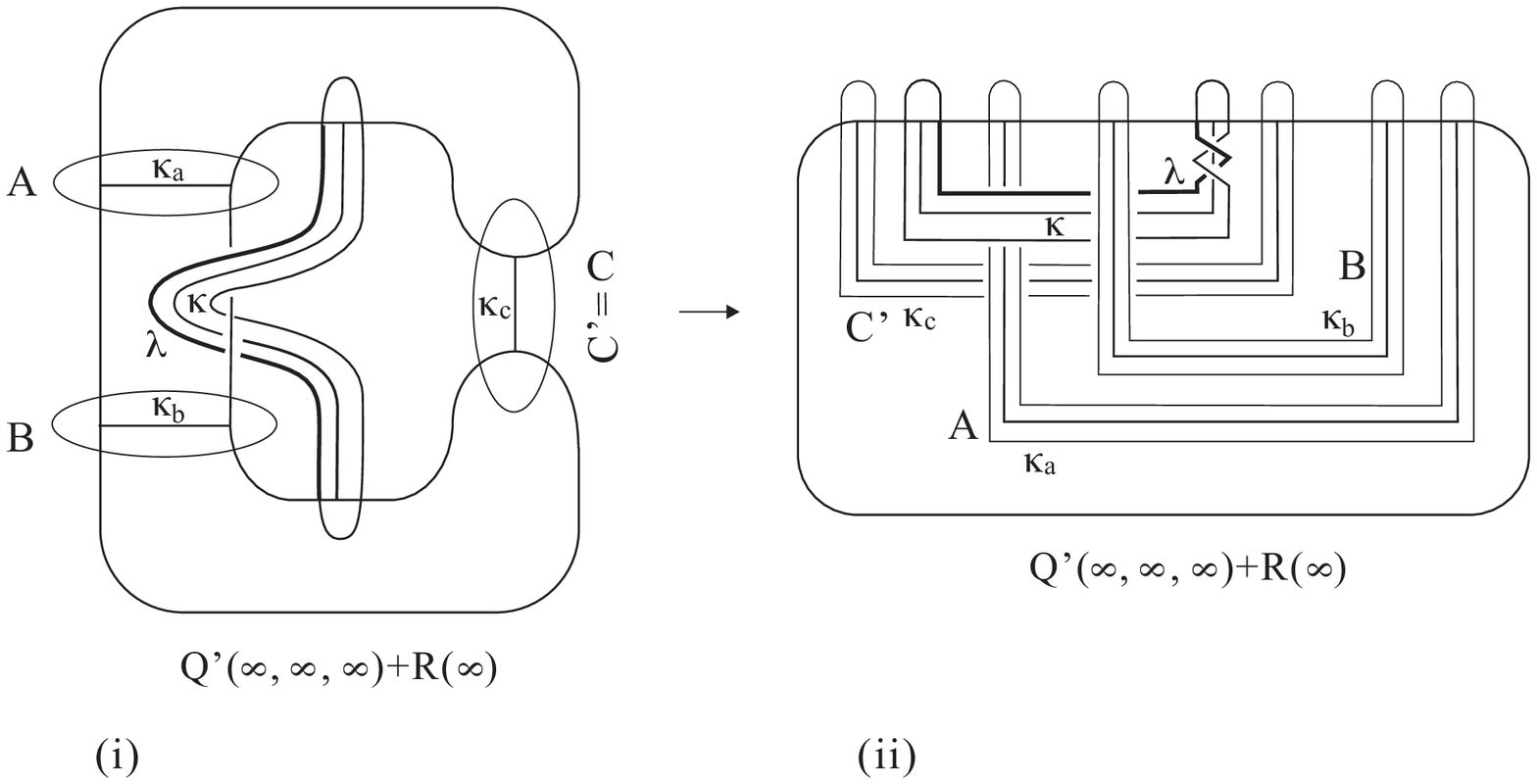}
\caption{}
\label{KABCcovering}
\end{center}
\end{figure}

Following \cite{EM2},
we obtain a surgery description of
the covering knot $K(A, B, C)$.
By an ambient isotopy of $S^3$ we move
$Q(\infty, \infty, 0)+R(\infty)$ to the position
in Figure~\ref{KABCcovering}(i);
during this isotopy the tangles $A, B, C$ are fixed.
We let $C'$ be the rational tangle obtained from $C$
by $(-\frac{\pi}{2})$--rotation about a line perpendicular to
the projection plane,
and define the tangle $Q'(A, B, C')$ to be $Q(A, B, C)$.
Note that
$Q(\frac{\alpha_1}{\beta_1},  \frac{\alpha_2}{\beta_2}, \frac{\alpha_3}{\beta_3}) 
= Q'(\frac{\alpha_1}{\beta_1}, \frac{\alpha_2}{\beta_2}, -\frac{\beta_3}{\alpha_3})$;
in particular, $Q(\infty, \infty, 0) = Q'(\infty, \infty, \infty)$.
After a further isotopy, 
we obtain Figure~\ref{KABCcovering}(ii).  
Note that $\kappa_a, \kappa_b, \kappa_c, \kappa$
in Figure~\ref{KABCcovering}
are spanning arcs of the $\infty$--tangles $A, B, C', R(\infty)$,
respectively.

\begin{figure}[htbp]
\begin{center}
\includegraphics[width=0.8\linewidth]{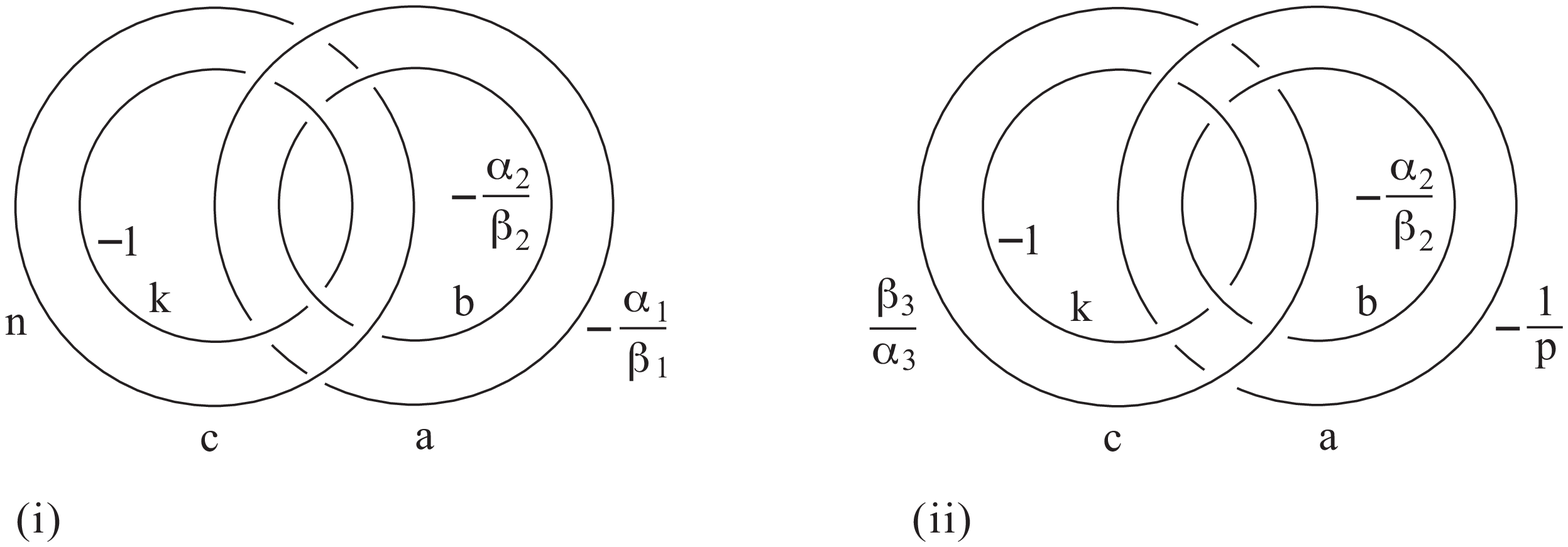}
\caption{}
\label{KABCcovering2}
\end{center}
\end{figure}

In the two--fold branched cover of $S^3$
along the trivial knot $Q'(\infty, \infty, \infty)+R(\infty)$,
we denote the preimages of $\kappa_a, \kappa_b, \kappa_c, \kappa$
by $a, b, c, k$, respectively.
We see from Figure~\ref{KABCcovering}(ii) that
the 4--component link $a\cup b\cup c\cup k$
is as illustrated in Figure~\ref{KABCcovering2}.
Furthermore,
the preimages of the latitudes of the $\infty$--tangles
$A, B, C'$ are preferred longitudes of $a, b, c$, respectively. 
Note that each of $a, b, c, k$ is
the covering knot of some trivializable tangle.
For example,
$a$ is the covering knot of
the tangle $Q'(\quad,\infty, \infty)+R(\infty)$ and
$k$ is the covering knot of $Q'(\infty, \infty, \infty)$.
Therefore,
$\frac{\alpha_1}{\beta_1}$--, $\frac{\alpha_2}{\beta_2}$--,
$(-\frac{\beta_3}{\alpha_3})$--untangle surgeries along 
$\kappa_a$, $\kappa_b$, $\kappa_c$ correspond to 
$(-\frac{\alpha_1}{\beta_1})$--, $(-\frac{\alpha_2}{\beta_2})$--,
$\frac{\beta_3}{\alpha_3}$--surgeries on 
$a$, $b$, $c$ in Figure~\ref{KABCcovering2}. 
We thus obtain the covering knot 
$K(\frac{\alpha_1}{\beta_1}, \frac{\alpha_2}{\beta_2}, \frac{1}{n})$
(resp.\
$K(\frac{1}{p}, \frac{\alpha_2}{\beta_2}, \frac{\alpha_3}{\beta_3})$) from $k$ in Figure~\ref{KABCcovering2}(i) (resp.\ (ii))
after the surgeries on $a, b, c$ given in
Figure~\ref{KABCcovering2}(i) (resp.\ (ii))\cite[Proposition 6.3]{EM2}.  
Note that the preimage of the latitude $\lambda$ of $R(\infty)$
in Figure~\ref{KABCcovering}(ii) 
gives the $(-2)$--framing of $k$.
Hence,
 $m$--untangle surgery along $\kappa$ corresponds
to $-(m+2)$--surgery on $k$; 
in particular, $(-1)$--untangle surgery on $\kappa$ corresponds to 
$(-1)$--surgery on $k$. 
We thus have a ``surgery description" of
the Seifert surgery $(K(A, B, C), \gamma)$ as follows.
For a link $k_1 \cup \cdots \cup k_n$ in $S^3$ and 
$r_i$ a slope on $\partial N(k_i)$, 
$(k_1, \dots, k_n;\ r_1, \dots, r_n)$ denotes
an $n$--tuple of $r_i$--surgeries
on $k_i$ $(1 \le i \le n)$.

\begin{proposition}
\label{surgery description KABC}
\begin{enumerate}
\item
The Seifert surgery 
$(K(\frac{\alpha_1}{\beta_1}, \frac{\alpha_2}{\beta_2}, \frac{1}{n}), \ \gamma)$ 
is obtained from $(k, -1)$
by the triple of surgeries 
$(a, b, c;\,
-\frac{\alpha_1}{\beta_1}, -\frac{\alpha_2}{\beta_2}, n)$
as in Figure~\ref{KABCcovering2}(i).

\item
The Seifert surgery 
$(K(\frac{1}{p}, \frac{\alpha_2}{\beta_2}, \frac{\alpha_3}{\beta_3}), \ \gamma)$  
is obtained from $(k, -1)$ by the triple of surgeries 
$(a, b, c;\,
-\frac{1}{p}, -\frac{\alpha_2}{\beta_2}, \frac{\beta_3}{\alpha_3})$
as in Figure~\ref{KABCcovering2}(ii).

\end{enumerate}
\end{proposition}

\begin{figure}[htbp]
\begin{center}
\includegraphics[width=0.8\linewidth]{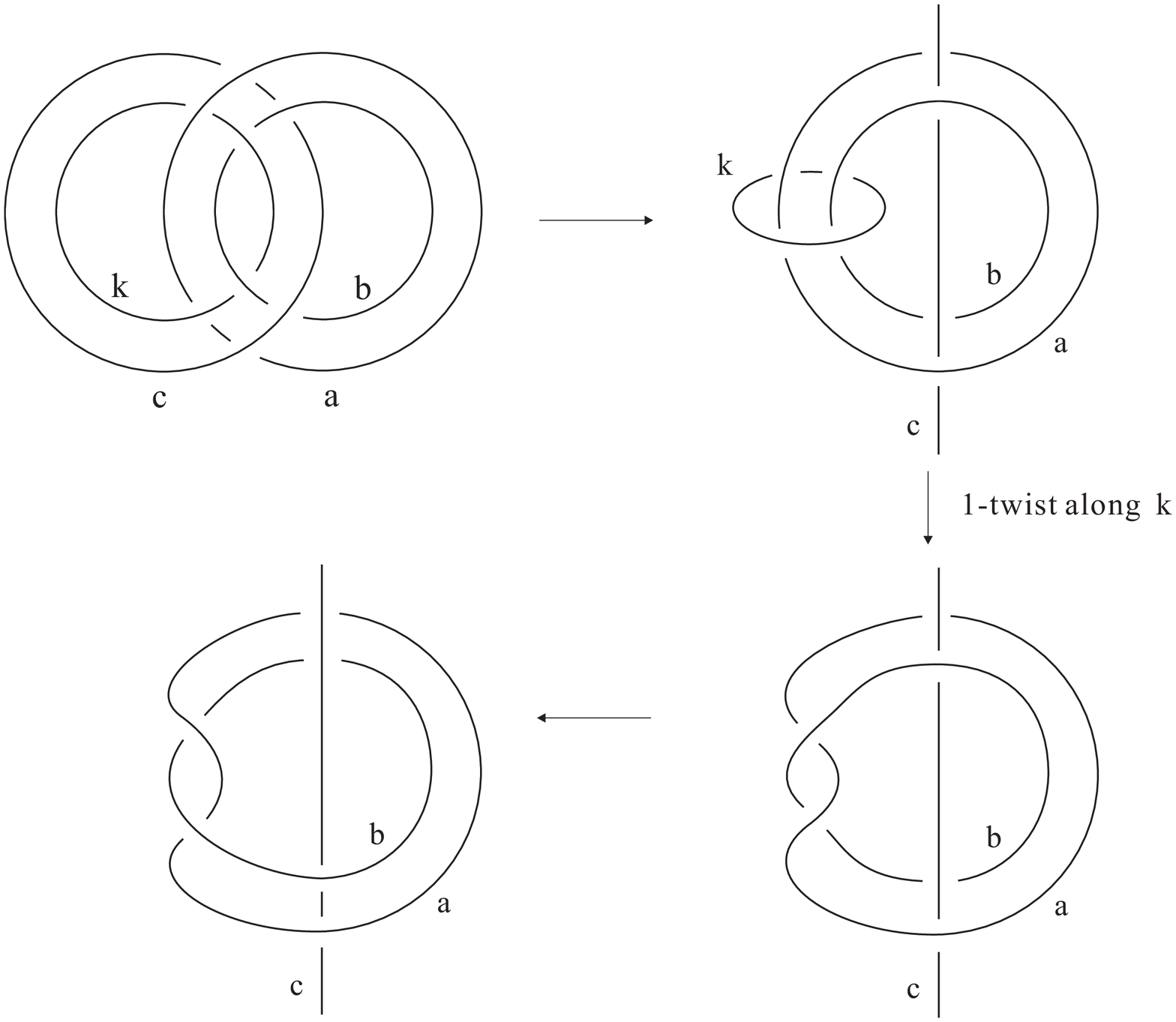}
\caption{$a, b, c$ are fibers in $k(-1) \cong S^3$ simultaneously.}
\label{tripleseiferters}
\end{center}
\end{figure}

Figure~\ref{tripleseiferters} shows that
$1$--twist along $k$ converts 
$a \cup b \cup c$ to a union of fibers in a Hopf fibration of $k(-1) \cong S^3$. 
We thus have the following.

\begin{proposition}
\label{seiferters for KABCi}
The knots $a, b, c$ are seiferters
for the Seifert surgery $(k, -1)$; 
more precisely, 
they become fibers in a Hopf fibration of $k(-1) \cong S^3$,  simultaneously. 
\end{proposition}

\medskip
We show that 
the Seifert surgery $(K(A, B, C), \gamma)$ is obtained
from $(k, -1)$ after a sequence of twistings along 
seiferters $a, b, c$.

\begin{lemma}[{\cite[Claim 6.10]{DMM1}}]
\label{splitting}
Let $a \cup b$ be a Hopf link. 
If a pair of surgeries 
$(a, b;\, \frac{x}{y}, \frac{s}{t})$ satisfies $| xs - yt | = 1$,
then it can be realized by a finite 
sequence of alternate twistings: 
$$(a,\,  m_1\mathrm{-twist}) \to (b,\, n_1\mathrm{-twist}) \to \cdots  \to (a,\, m_p\mathrm{-twist}) 
\to (b,\, n_p\mathrm{-twist})$$
for some integers $m_i, n_i$, 
where only $n_p$ can be $0$. 
\end{lemma}

\textsc{Proof.}
We first note that 
after 
arbitrary twists along $a$ or $b$, 
$a\cup b$ remains a Hopf link. 
Hence, we can keep twisting along $a$ and $b$ alternately.
Note also that
$m$--twist along $a$ changes $(a, b;\, \frac{x}{y}, \frac{s}{t})$
to $(a, b;\, \frac{x}{y+mx}, \frac{s+mt}{t})$ and 
$n$--twist along $b$ changes $(a, b;\, \frac{x}{y}, \frac{s}{t})$
to $(a, b;\, \frac{x+ny}{y},\frac{s}{t+ns})$; 
see \cite[9.H]{Ro}. 
The Euclidean algorithm changes $(x, y)$ to $(1, 0)$,
so that alternate twistings 
along $a$ and $b$ convert 
the pair of surgery $(a, b;\, \frac{x}{y}, \frac{s}{t})$
to $(a, b;\, \frac{1}{0}, \frac{s'}{t'})$. 
Since $|1\cdot s' - 0 \cdot t'| = 1$, 
we may assume $s' = 1$. 
Thus applying $(-t')$--twist along $b$, 
we have $(a, b;\, \frac{1}{0}, \frac{1}{0})$. 
This implies the required result. 
\hspace*{\fill} $\square$(Lemma~\ref{splitting})

\medskip

\begin{remark}
\label{splitting_remark}
In Lemma~\ref{splitting}, 
a finite sequence of twistings realizing a pair of surgeries 
$(a, b;\, \frac{x}{y}, \frac{s}{t})$ is not unique. 
For instance, 
$(a, b;\, \frac{x}{y}, \frac{s}{t})$ is realized also by 
a finite sequence of alternate twistings: 
$$(b,\, m_1\mathrm{-twist}) \to (a,\, n_1\mathrm{-twist}) \to \cdots \ \to (b,\, m_p\mathrm{-twist}) 
\to (a,\, n_p\mathrm{-twist})$$
for some integers $m_i, n_i$, 
where only $n_p$ can be $0$. 
\end{remark}

\begin{figure}[htbp]
\begin{center}
\includegraphics[width=1.0\linewidth]{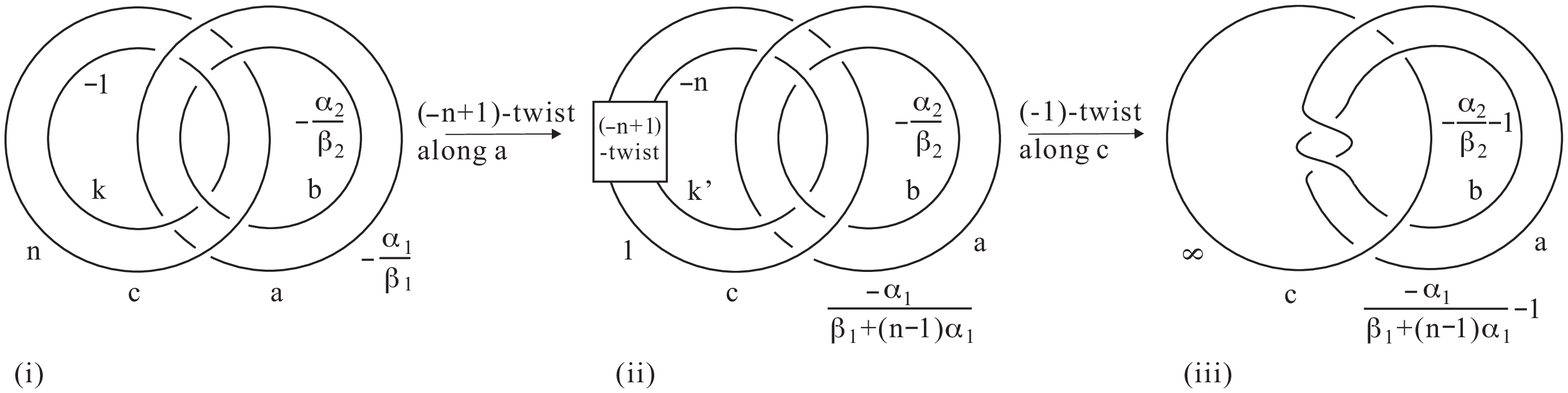}
\caption{} 
\label{surgeryDescriptionKn}
\end{center}
\end{figure}

First we consider the Seifert surgery
$(K, \gamma) =
(K(\frac{\alpha_1}{\beta_1}, \frac{\alpha_2}{\beta_2}, \frac{1}{n}), \gamma)$ 
satisfying Proposition~~\ref{TABC}(1)(i). 
By Proposition~\ref{surgery description KABC}(1)
the Seifert surgery $(K, \gamma)$ is obtained from
the Seifert surgery $(k, -1)$ 
after the surgeries on $a, b, c$ given in Figure~\ref{surgeryDescriptionKn}(i). 
Figure~\ref{surgeryDescriptionKn}(ii) gives
the surgery description of $(K, \gamma)$ 
after $(-n +1)$--twist along the seiferter $a$,
where the images of $a, b, c, k$ are denoted by
$a, b, c, k'$, respectively.
We note that $a, b, c$ remain seiferters for the Seifert surgery
$(k', -n)$.
Figure~\ref{surgeryDescriptionKn}(ii) indicates that
$(K, \gamma)$ is obtained from $(k', -n)$
by the triple of surgeries
$(a, b, c;\, \frac{-\alpha_1}{\beta_1 + (n-1)\alpha_1}, -\frac{\alpha_2}{\beta_2}, 1)$. 
In Figure~\ref{surgeryDescriptionKn}(ii), 
$k'$ is a $(1, -n+1)$ cable of the solid torus $S^3 -\mathrm{int}N(c)$,
and $a$ is the core of the solid torus.
Hence, performing the $1$--surgery (i.e.\ $(-1)$--twist) on $c$
given in Figure~\ref{surgeryDescriptionKn}(ii) converts
$k' = T_{1, -n+1}$ to $T_{1-(-n+1), -n+1} = T_{n, -n+1}$,
and the surgery coefficient $-n$ on $k'$ to
$-n - (-n+1)^2 = n(1-n)-1$;
the surgery coefficients on $a$ and $b$ decrease by one.
It follows that $(K, \gamma)$ is obtained from
$(T_{n, -n+1}, n(1-n)-1)$ by the pair of surgeries
$(a, b;\,  \frac{-\alpha_1}{\beta_1 + (n-1)\alpha_1}-1, -\frac{\alpha_2}{\beta_2}-1)$.

\begin{figure}[htbp]
\begin{center}
\includegraphics[width=0.9\linewidth]{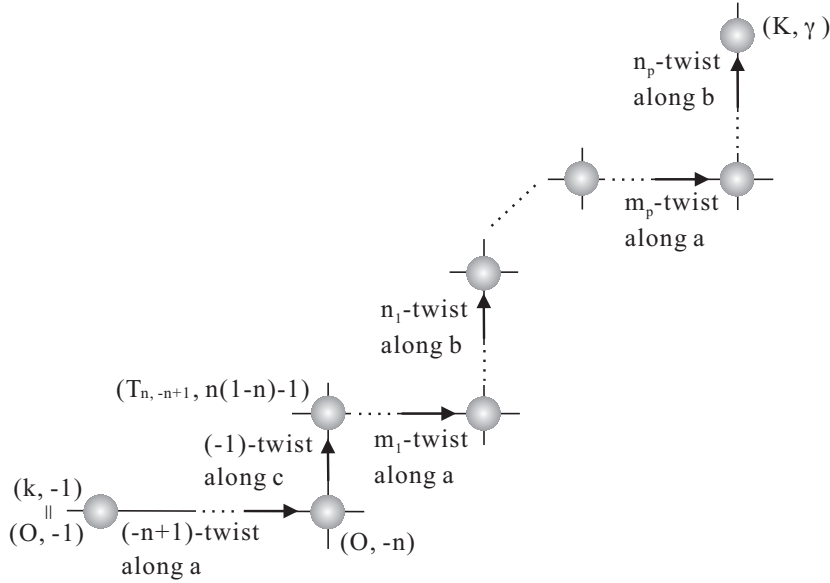}
\caption{$(K, \gamma) =
(K(\frac{\alpha_1}{\beta_1}, \frac{\alpha_2}{\beta_2}, 
\frac{1}{n}), \gamma)$.} 
\label{PortionKABn}
\end{center}
\end{figure}

Note that $\{a, b\}$ is a pair of seiferters
for $(T_{n, -n+1}, n(1-n)-1)$ and forms a Hopf link.
On the other hand,
we have $n\alpha_1\alpha_2 + \alpha_1 \beta_2 + \beta_1 \alpha_2 = \pm 1$ by Proposition~\ref{TABC}(1)(i), so that
$(a, b;\,  \frac{-\alpha_1}{\beta_1 + (n-1)\alpha_1}-1, -\frac{\alpha_2}{\beta_2}-1)$ 
satisfies the condition of Lemma~\ref{splitting}.
It follows that
$(K, \gamma)$ 
is obtained from $(T_{n, -n+1}, n(1-n)-1)$ by applying 
alternate twistings: 
$$(a,\,  m_1\mathrm{-twist}) \to (b,\,  n_1\mathrm{-twist}) \to \cdots \to (a,\,  m_p\mathrm{-twist}) 
\to (b,\,  n_p\mathrm{-twist})$$
for some integers $m_i, n_i$.  
We then obtain
Figure~\ref{PortionKABn} which gives an explicit path from 
$(T_{n, -n+1}, n(1-n)-1)$ (and also $(O, -1)$)
to $(K, \gamma)$ in the Seifert Surgery Network. 
\hspace*{\fill}
$\square$(Theorem~\ref{connected}
for $(K(\frac{\alpha_1}{\beta_1}, \frac{\alpha_2}{\beta_2}, \frac{1}{n}), \gamma)$ in $\mathcal{EM}\mathrm{III}$)

\begin{remark}
\label{another path}
As mentioned in Remark~\ref{splitting_remark}, 
a finite sequence of twistings along $a$ and $b$ realizing 
$(a, b; \frac{x}{y}, \frac{s}{t})$ with $| xs - yt | = 1$
is not unique, and thus
there are other paths from $(T_{n, -n+1}, n(1-n)-1)$ to 
$(K, \gamma)$. 
\end{remark}

\begin{figure}[htbp]
\begin{center}
\includegraphics[width=1.0\linewidth]{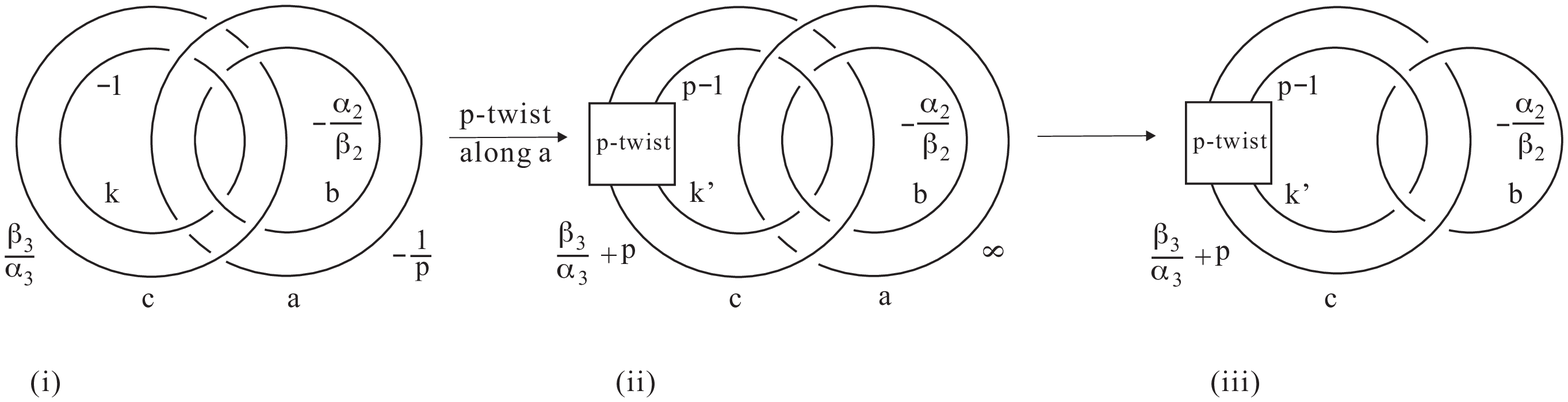}
\caption{}
\label{surgeryDescriptionKp}
\end{center}
\end{figure}

Finally we consider $(K, \gamma) =
(K(\frac{1}{p}, \frac{\alpha_2}{\beta_2}, \frac{\alpha_3}{\beta_3}), \gamma)$ 
satisfying Proposition~~\ref{TABC}(1)(ii). 
Performing the $(-\frac{1}{p})$--surgery
on $a$ in Figure~\ref{surgeryDescriptionKp}(i),
we obtain a new surgery description of $(K, \gamma)$ 
as in Figure~\ref{surgeryDescriptionKp}(ii),  
in which $k'$, the image of $k$,
remains unknotted in $S^3$ and 
the surgery coefficient of $k'$ is $p-1$.   
Note that $\{b, c\}$ is a pair of seiferters
for $(k', p-1)$.
It follows that
$(K, \gamma)$ is obtained from $(k', p-1)$ after a pair of surgeries
on the Hopf link $b\cup c$ given in Figure~\ref{surgeryDescriptionKp}(iii).

\begin{figure}[htbp]
\begin{center}
\includegraphics[width=0.77\linewidth]{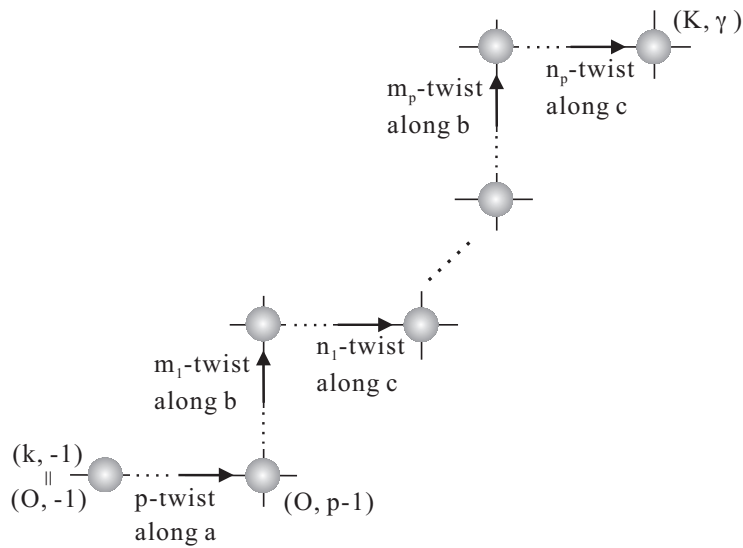}
\caption{$(K, \gamma) = 
(K(\frac{1}{p}, \frac{\alpha_2}{\beta_2}, \frac{\alpha_3}{\beta_3}), \gamma)$.}
\label{portionKpBC}
\end{center}
\end{figure}

As in the previous case, 
since $p\alpha_2\alpha_3 + \alpha_2\beta_3 + \beta_2\alpha_3
= \pm 1$ by Proposition~\ref{TABC}(1)(ii), 
the pair of surgeries  
$(b, c;\, -\frac{\alpha_2}{\beta_2}, \frac{\beta_3}{\alpha_3} + p)$
in Figure~\ref{surgeryDescriptionKp}(iii)
satisfies the condition of Lemma~\ref{splitting}. 
Thus $(K, \gamma)$ is obtained from $(O, p-1)$ after
alternate twistings along the seiferters $b, c$: 
$$(b,\, m_1\mathrm{-twist}) \to (c,\, n_1\mathrm{-twist}) \to \cdots  \to (b,\, m_p\mathrm{-twist})  
\to (c,\, n_p\mathrm{-twist})$$
for some integers $m_i, n_i$. 
See Figure~\ref{portionKpBC}. 
\hspace*{\fill} $\square$(Theorem~\ref{connected}
for $(K(\frac{1}{p}, \frac{\alpha_2}{\beta_2}, \frac{\alpha_3}{\beta_3}), \gamma)$ in $\mathcal{EM}\mathrm{III}$) 

As mentioned in Remark~\ref{another path} for 
$(K(\frac{\alpha_1}{\beta_1}, \frac{\alpha_2}{\beta_2}, 
\frac{1}{n}), \gamma)$, 
there is yet another path from $(O, p-1)$ to 
$(K(\frac{1}{p}, \frac{\alpha_2}{\beta_2}, \frac{\alpha_3}{\beta_3}), \gamma)$. 

\end{document}